\documentclass[reqno,twoside]{amsart}
\usepackage{amsmath}
\usepackage{amsfonts}
\usepackage{amssymb}
\usepackage{enumerate}

\usepackage{a4wide,amsmath,amssymb,latexsym,amsthm}

\usepackage{mathrsfs,mathtools,epic,bm}
\usepackage{hyperref}
\hypersetup{colorlinks=true,linkcolor=blue,filecolor=mangeta,urlcolor= cyan}

 \usepackage[pagewise]{lineno}

\newtheorem{theorem}{Theorem}[section]
\newtheorem{proposition}[theorem]{Proposition}
\newtheorem{lemma}[theorem]{Lemma}

\theoremstyle{definition}
\newtheorem{definition}[theorem]{Definition}

\newtheorem{remark}[theorem]{Remark}

\newtheorem{assumption}[theorem]{Assumption}

\newcommand\norm[1]{\left\lVert #1\right\rVert}
\newcommand\abs[1]{\left\lvert#1\right\rvert}

\numberwithin{equation}{section}

\def \dis {\displaystyle}

\def \R {\mathbb{R}}

\def \Ha {{H^{\alpha}_{a^+}(a, b)}}
\def \Hi {{H^{\alpha}_{a^+}(a, b_i)}}
\def \Hb {{H^{\alpha}_{b^-}(a, b)}}

\def\hy {\hat{y}}
\def\hu {\hat{u}}
\def\hp {\hat{p}}

\def \Dc {{\mathcal D}}
\def \Dr {{\mathbb D}}
\def \D {{\mathbb D}}

\def \Da { \Dr_{a^+}^\alpha}

\keywords{Riemann-Liouville fractional derivative, Caputo fractional derivative, fractional integral, Sturm-Liouville equations, initial-boundary value problems, optimal control, optimality system, optimaltiy conditions}

\begin{document}
\title[Fractional Sturm-Liouville equations on a star graph]{Optimal control problems of parabolic fractional Sturm-Liouville equations in a star graph }

\author{G\"unter Leugering}
\address{G\"unter Leugering, Department of Mathematik, Friedrich-Alexander-Universit\"at Erlangen-N\"urnberg, Cauerstr. 11 (03.322), 91058 Erlangen, Germany.}
\email{guenter.leugering@fau.de}

\author{Gis\`{e}le Mophou}
\address{Gis\`{e}le Mophou, Laboratoire L.A.M.I.A., D\'{e}partement de Math\'{e}matiques et Informatique, Universit\'{e} des Antilles, Campus Fouillole, 97159 Pointe-\`{a}-Pitre,(FWI), Guadeloupe -
Laboratoire  MAINEGE, Universit\'e Ouaga 3S, 06 BP 10347 Ouagadougou 06, Burkina Faso.}
\email[Mophou]{gisele.mophou@univ-antilles.fr}

\author{Maryse Moutamal}
\address{Maryse Moutamal, University of Buea, Department of Mathematics, Buea, Cameroon.}
\email{maryse.moutamal@aims-cameroon.org}

\author{Mahamadi Warma}
\address{Mahamadi Warma, Department of Mathematical Sciences and the Center for Mathematics and Artificial Intelligence (CMAI), George Mason University,  Fairfax, VA 22030, USA.}
\email{mwarma@gmu.edu}


\thanks{The third author is supported by the Deutscher Akademischer Austausch Dienst/German Academic Exchange Service (DAAD). The fourth author is partially supported by the AFOSR under Award NO:  FA9550-18-1-0242 and by the US Army Research Office (ARO) under Award NO: W911NF-20-1-0115.}


\begin{abstract}	
In the present paper we deal with parabolic fractional initial-boundary value problems of Sturm–Liouville type  in an interval and  in a general star graph. We first give several existence, uniqueness and regularity results of weak and very-weak solutions.
We prove the existence and uniqueness of solutions to a quadratic boundary optimal control problem and provide a characterization of the optimal contol via the  Euler–Lagrange first order optimality conditions.  We then investigate the analogous problems for a fractional
Sturm–Liouville problem in a general star graph with mixed Dirichlet and Neumann boundary controls. The existence and uniqueness of minimizers, and the characterization of the first order optimality conditions are obtained in a general star graph by using the method of Lagrange multipliers.
\end{abstract}

\maketitle

\section{Introduction and problem setting}

\begin{figure}[h]
	\centering
	\includegraphics[width=0.8\linewidth]{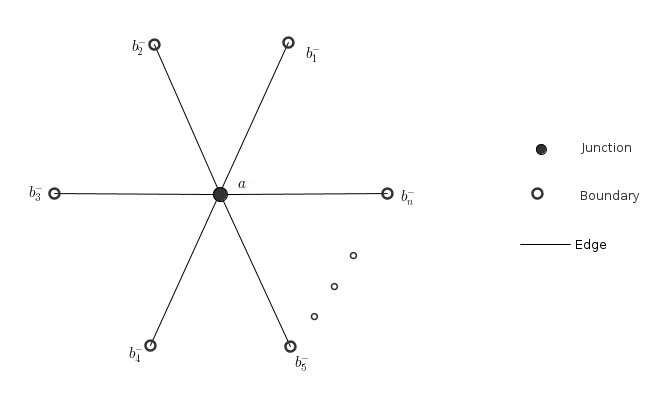}
	\caption{A sketch of a star graph with n edges}
	\label{fig:stargraph}
\end{figure}

The  main concern of the present paper is to study fractional optimal control problems on a link to a star graph (see Figure \ref{fig:stargraph}) by considering the following minimization problem:
\begin{align}\label{OP2}
\underset{v\in {\mathbb{U}_{ad}}}{\min}\left(\frac{1}{2}\sum_{i=1}^{n}\int_{Q_i} \abs{y^i-y_d^i}^2\,dxdt+\frac{1}{2}\sum_{i=2}^{m}\int_0^T|u_i|^{2} \,dt +\frac{1}{2}\sum_{i=m+1}^{n}\int_0^T|v_i|^{2} \,dt\right),
\end{align}
subject to the constraints that $y=(y^i)_i$ ($i=1,2,\ldots,n$,  where $n\in\mathbb N$ is fixed) satisfies the following parabolic system involving a fractional  Sturm-Liouville operator on a star graph:
\begin{equation}\label{ST2}
	\left\{
\begin{array}{lllllllllllllllll}
\displaystyle	y_t^i+\mathcal{D}_{b_i^-}^\alpha(\beta^i\mathbb{D}_{a^+}^\alpha y^i)+q^iy^i&=&f^i&\text{in}& Q_i:=(a,b_i)\times (0,T),\,i=1,\dots, n,\\
\displaystyle	(I_{a^+}^{1-\alpha}y^i)(a,\cdot)-(I_{a^+}^{1-\alpha}y^j)(a,\cdot)&=&0&\text{in}&(0,T),~i\neq j=1,\dots, n,\\
\displaystyle	\sum_{i=1}^n(\beta^i\mathbb{D}_{a^+}^\alpha y^i)(a,\cdot)&=&0&\text{in}& (0,T), \\
\displaystyle	(I_{a^+}^{1-\alpha} y^1)(b_1^-,\cdot)&=&0&\text{in}&(0,T), \\
\displaystyle	(I_{a^+}^{1-\alpha} y^i)(b_i^-,\cdot)&=&u_i&\text{in}& (0,T), \; i=2,\dots, m\\
\displaystyle	(\beta^i\mathbb{D}_{a^+}^\alpha y^i)(b_i^-,\cdot)&=&v_i&\text{in}&(0,T), \; i=m+1,\dots,n,\\
\displaystyle	y^i(\cdot,0) &=&0&\text{in}& (a,b_i),~~i=1,\dots,n,
	\end{array}
	\right.
\end{equation}
where  $y_d=(y_d^i)_i\in L^2\left((0,T);\left(L^2(a,b_i)\right)^n\right)$, $\mathbb{U}_{ad}$ is a closed and convex subset of $ \left(L^2(0,T)\right)^{n-1}$,  and $2\leq m\leq  n$ is a natural number.
Here, $T>0,$ $a, b_i\in\mathbb R$ with $0\leq a<b_i$, $\Dr_{a^+}^\alpha $,  $\Dc_{b_i^-}^\alpha,\, i=1,\dots,n,$ stand  for  the left Riemann-Liouville,  and the right Caputo fractional derivatives of order $\alpha \in (0,1]$, respectively, 
and $I^\alpha_{a^+}$ ($0<\alpha\leq 1$) is the Riemann-Liouville fractional integral of order $\alpha\in (0,1]$.  We refer to Section \ref{prelim} for the precise definition.  
The real valued functions $\beta^i\in C([a,b_i])$ and $q^i\in L^\infty(a,b_i),\, i=1,\dots,n$,  satisfy suitable conditions (see Assumption \ref{assump2} below), $f^i$ belongs to $L^2(Q_i), _, i=1,\dots,n$,   the controls $u_i\in L^2(0,T)$,  $i=2,\ldots m$, and $v_i\in L^2(0,T), i=m+1,\dots,n$.  Some controls can be equal to zero.  In addition,  if $c=a$ or $c=b_i^-$, then $(I_{a^+}^{1-\alpha} y)(c,\cdot)=\lim_{x\to c}(I_{a^+}^{1-\alpha} y)(x,\cdot)$,  and $(\beta\mathbb{D}_{a^+}^\alpha y)(c,\cdot)=\lim_{x\to c}(\beta\mathbb{D}_{a^+}^\alpha y)(x,\cdot)$ for smooth functions $y$, otherwise it is understood in the weak sense as in the formulation of our notion of weak solutions (see Definitions \ref{weaksolutionnh} and \ref{weaksolutionnonhomostar}).

The setting indicates that we are looking at a star graph, rooted at $b^-_1$, i.e., where we have a fixed Dirichlet-type boundary condition and controls via fractional Dirichlet and Neumann boundary conditions.
After proving some well-posedness results (existence and uniqueness of weak and very-weak solutions) of the system \eqref{ST2} in the general star graph,  we show the existence and uniqueness of minimizers to the optimal control problem \eqref{OP2}-\eqref{ST2},  and give  the associated optimality conditions by using the method of Lagrange multipliers.

Notice that our boundary conditions are given in terms of the Riemann-Liouville fractional integral.
On the contrary,  boundary (initial) conditions for the Caputo derivatives  are expressed
in terms of boundary (initial)  values of integer order derivatives. This allows for a
numerical treatment of initial value problems for differential equations of non
integer order independently of the chosen definition of the fractional derivative.
For this reason, many authors either resort to Caputo derivatives, or
use the Riemann-Liouville derivatives but avoid the problem of boundary (initial) values
of fractional integrals by treating only the case of zero boundary (initial) conditions.  
The interesting paper \cite{Pod2} has provided a series of examples from the field of viscoelasticity which demonstrates
that it is possible to attribute physical meaning to boundary (initial) conditions
expressed in terms of Riemann-Liouville fractional integrals (as in \eqref{ST2} and \eqref{ST1} below),
and that it is possible to obtain boundary (initial) values for such initial conditions by appropriate measurements or observations. The mentioned examples include: The Spring-pot model, which is a linear viscoelastic element whose behavior is intermediate between that of an elastic element 
and a viscous element; a stress relaxation or a general deformation; and an impulse response.  For more details we refer to \cite{Pod2} and the references therein.

In order to tackle the above problem we need some preparation.  We first consider the following optimal control problem:
\begin{align}\label{OP1}
\dis \min_{v\in {\mathcal{U}_{ad}}}\left(\frac 12\int_Q \abs{y-y_d}^2\,dxdt+\frac{N}{2}\int_0^T|v|^2 \,dt\right),
\end{align}
subject to the constraint that the state $y=y(v)$ satisfies  the following initial-boundary value fractional Sturm–Liouville parabolic equation:
\begin{equation}\label{ST1}
\left\{\begin{array}{lllllll}
\displaystyle  y_t+\dis  \Dc_{b^-}^\alpha\,(\beta\, \Dr_{a^+}^\alpha y)+q\,y&=&f&\hbox{ in } Q:=(a,b)\times(0,T),\\
\displaystyle (I_{a^+}^{1-\alpha} y)(a,\cdot)&=& 0&\hbox{ in } (0,T),\\
\displaystyle (\beta\Dr_{a^+}^{\alpha}y)(b^-,\cdot)&=&v&\hbox{ in } (0,T),\\
\displaystyle y(\cdot,0)&=&y^0&\hbox{ in } (a,b).
\end{array}
\right.
\end{equation}
Here, $y_d\in L^2(Q),$  $N>0$ is a real number and $\mathcal{U}_{ad}$ is a closed and convex subset of $L^2(0,T)$.
The real number  $T>0$, $a,b\in\mathbb R$ with $0\leq a<b$,  $\Dr_{a^+}^\alpha $, $\Dc_{b^-}^\alpha$, and $I^\alpha_{a^+}$ are as above.  
The real valued functions $\beta\in C([a,b])$ and $q\in L^\infty(a,b)$ satisfy suitable conditions (see Assumption \ref{asump-beta-q}),  $f\in L^2(Q)$,  and the control $v\in L^2(0,T)$.   Here also,  $(I_{a^+}^{1-\alpha} y)(a,\cdot)$ and $(\beta\Dr_{a^+}^{\alpha}y)(b^-,\cdot)$ are interpreted as above. 
After proving several existence, uniqueness and regularity results of the state equation \eqref{ST1} and the associated dual system,  we show the existence and uniqueness of minimizers of the optimal control problem \eqref{OP1}-\eqref{ST1}, and characterize the associated first order optimality conditions by using the classical Euler-Lagrange first order optimality conditions.

The results obtained in the present work generalize to the parabolic setting,  the ones contained in \cite{mophou2020} for fractional elliptic Sturm-Liouville problems.

Classical Sturm–Liouville theory is the study of second-order linear ordinary differential equations of the form:
\begin{equation*}
	(\beta y')'+qy= \lambda \omega y
\end{equation*}
where $y$ is the unknown, that is,  a physical quantity, $\lambda$ a parameter and $\beta, q,\omega$ are suitable functions.  The function $\omega(x)$ is called the weight or density function.
With appropriate boundary conditions,  $\lambda$ and $y$ appear as  eigenvalue and eigenfunction, respectively,  of the adjoint operator. There is a wide literature related to this kind of equations. We refer for instance  to \cite{Zettl} and the references therein.

It has been amply shown in the past that many phenomena which occur in various fields in science as well as in engineering can be more accurately described by means of fractional order derivatives than integer order derivatives.  Many researchers have then focused their attention on fractional Sturm–Liouville problems which are obtained by replacing the ordinary derivatives with fractional derivatives. For instance,  \cite{Zayer}  studied a fractional Sturm-Liouville eigenvalue problem involving a right-sided Riemann-Liouville fractional derivative,  and a left-sided Caputo fractional derivative of the same order. It has been shown that analytical solutions are non-polynomial functions and are orthogonal with respect to the weighted function associated to the problem.  The authors in \cite{klimek2013}  proved the orthogonality of solutions  to  fractional Sturm-Liouville eigenvalue problems involving first, a left-sided Riemann-Liouville fractional derivative, a right-sided Caputo fractional derivative,  and second a right-sided Riemann-Liouville fractional derivative,  and a left-sided Caputo fractional derivative of the same order. Using a variational approach,  the authors in \cite{klimek2014} showed the existence of a countable set of orthogonal solutions and eigenvalues to a fractional Sturm-Liouville eigenvalue problem involving left-sided and right-sided Caputo fractional derivative of the same order.  The paper \cite{rivero2013} provided an approach to the fractional version of the Sturm-Liouville problem, by using different fractional operators that coincide to the ordinary operator when the fractional parameter becomes an integer.  Moreover, for each fractional operator, some basic properties of the Sturm-Liouville theory have been investigated. We also refer to the work \cite{hassan2015}, where a fractional Sturm-Liouville eigenvalue problem is studied in an unbounded domain,  and to the work \cite{qasem} for the numerical approximation of such problems.  The work \cite{darius} used Riemann–Liouville fractional derivatives to introduce and characterize fractional Sobolev spaces which, in turn, have been utilized, via a variational approach, for the existence and uniqueness  of solutions to a boundary value problem associated to a Sturm-Liouville type equation involving left-sided and right-sided  Riemann-Liouville fractional derivatives  of the same order.

Fractional order operators are typical examples that model anomalous diffusion.  It is nowadays known that anomalous diffusion  appear in many branches of science and engineering: for example in fluid pressure transients in porous media, thermal processes such as heat conduction in materials, or transport of chemicals and pollutants in the environment \cite{Selvadurai}.
The processes referring to the mentioned phenomena could be effectively modeled by fractional differential equations.  We notice that fractional derivatives are non-local operators. We refer to the monograph \cite{pod} for more details about fractional operators. Such characteristics allow to incorporate different types of information and to use  fractional derivatives to model systems with long-range interactions in space and time (memory effect).
This close link between diffusion systems and fractional derivatives gives rise to the increasing number of papers on the subject of fractional diffusion equations (see e.g., \cite{klimek2016,JNakagawa} and their references). In \cite{Chen},  the authors developed a theory of strong solutions of space–time fractional diffusion equations on bounded domains, as well as probabilistic representations of these solutions, which are useful for particle tracking codes.  The case of semilinear space-time fractional diffusion equations and systems has been thoroughly studied in the monograph \cite{Ga-Wa-2021}.  Super diffusive fractional evolution equations have been studied in \cite{AGKW-NA} and their references.
In \cite{klimek2016},  the authors used space-time fractional diffusion equations on finite domains to model anomalous diffusion behavior with large particle jumps combined with long waiting times. In their work, using the method of separating variables and applying theorems ensuring the existence of solutions to the fractional Sturm–Liouville problem,  they solved several types of fractional diffusion equations.

In the area of optimal control of fractional differential equations, there are also some interesting results: Agrawal detained the first record of the formulation of the fractional optimal control problem. He presented in \cite{agrawal2004general}  a general formulation and proposed a numerical method to solve such problems. In his work, the fractional derivative was defined in the Riemann–Liouville sense,  and the formulation was obtained through the fractional variation principle and the Lagrange multiplier technique.  The author in \cite{Mophou2011} applied the classical control theory to a fractional diffusion equation involving a Riemann–Liouville fractional derivative in a bounded domain by interpreting the Euler–Lagrange first-order optimality condition with an adjoint problem defined through a right fractional Caputo derivative. The author obtained an optimality system for the optimal control.  In \cite{dorville2011optimal},  a non-homogeneous Dirichlet boundary fractional diffusion equation in a bounded domain has been studied.  By transposition,  the existence and uniqueness of solutions of the boundary fractional diffusion equation has been shown,  and then under some appropriate assumptions on the closed convex set of the admissible controls,  a decoupling of the optimality system has been obtained.

Differential equations on graphs have many applications in mathematics, physics,  and engineering (see e.g., the monograph \cite{BK-2013} and the references therein). Several phenomena such as the flow on nets of gas pipeline \cite{Ste}, propagation of water waves in channel networks (the well-known Burgers equation)  lead to partial differential equations on graphs. Optimal control problems for classical partial differential equations have been investigated by many authors. For example the networks of beams and strings and their control properties have been studied by Lagnese el al. \cite{LG-2004}. The survey paper by D\'ager and Zuazua \cite{DZ-2006} is an excellent reference where several $1$-D optimal control problems on graphs have been provided.

The implementation of differential equations to the network domains, in particular in biology,  engineering,  cosmology,  leads Lumer \cite{Lumer} to the notion of evolution problems on "ramified spaces". Ever since, work on ordinary and partial differential equations on metric graphs or networks has greatly evolved (see e.g. \cite{Mehmeti, gunter1,dager, gunter2,schimidt,Vonbelow,gunter3} and their references). There are few results on  fractional differential equations on metric graphs or networks.  In \cite{mehra2019}, the authors studied  the existence and uniqueness of a nonlinear Caputo fractional boundary value problem on a star graph by means of a fixed point theory.  They obtained an equivalent system of fractional boundary value problems through  a proper transformation. Then they established existence and uniqueness results by a fixed point theory. This work was extended to the existence and uniqueness of solutions of a nonlinear fractional boundary value problem on a circular ring with an attached edge in \cite{mehra2021} by the same authors. The results are achieved by the Banach contraction principle and Krasnoselskii's fixed point theorem. Recently,  the paper \cite{mophou2020} considered  the optimal control of an  elliptic   problem involving a fractional  Sturm-Liouville operator in an interval,  and on a general star graph with mixed Dirichlet and Neumann boundary controls,  where the existence of solutions to a quadratic boundary optimal control problem and the associated optimality conditions have been proven.  In \cite{mehra2021a}, the authors investigated optimal control problems for nonlinear fractional order boundary value problems on a star graph,  where the fractional derivative is described in the Caputo sense. They derived the fractional optimality system by means of  the Lagrange multiplier method and  the Banach contraction principle.  Finally, in \cite{MML}, an optimal control problem for a time-fractional
diffusion process in a star graph has been considered, where the right-Caputo fractional derivative has been used with respect to the time. Existence and  uniqueness of solutions are provided for the system of equations as well as for the optimal control problem along with first order optimality conditions. Moreover,  a numerical simulation has been provided based on a finite difference discretization with respect to time and space.\\

The rest of the paper is structured as follows. In Section \ref{prelim}, we fix some notations,  give some definitions, introduce the function spaces needed to study  our problem and prove some preliminary results that will be used in the proofs of our main results.  In Section \ref{oneedge},  we show first that the homogeneous and non-homogeneous fractional Sturm-Liouville equations on a single edge have unique weak solutions.  The results are contained in Theorems \ref{Thm1} and \ref{theoexistnh}, respectively.
The regularity of solutions is also investigated.  We conclude this section by proving that the  quadratic optimal control problem associated to  the evolution equation involving a fractional Sturm-Liouville operator on one edge admits a unique optimal control,  and we give the associated optimality system that  characterizes this control (see Theorems \ref{theoopt} and \ref{theo-39}).  The same investigation is done for the evolution equation involving a fractional Sturm-Liouville operator on the considered graph in Section \ref{graph}.  In that case, we have introduced the notion of very weak-solutions and studied their existence and regularity.  The main results of this section are contained in Theorems \ref{theoexisthomostar},  \ref{theoexisthomostarnh1}, \ref{theo-48}, and Proposition \ref{regular}.

\section{Preliminaries}\label{prelim}

In this section, we introduce some notations,  give the function spaces needed to study our problems,  recall some known results and prove some intermediate results that are needed in the proofs of our main results. We start with fractional integrals and derivatives.

Let $[a,b]\subset\R$, $a\ge 0$ and $\rho : [a,b] \to \R$ be a given function.

\begin{definition}\label{defIRL}
The left,  and right Riemann–Liouville fractional integrals of order $\alpha\in(0,1]$ of $\rho\in L^1(a,b)$,  are defined, respectively,  by:
$$(I^{\alpha}_{a^+}\rho)(x):=\dis \frac{1}{\Gamma(\alpha)}\int_{a}^{x}(x-t)^{\alpha-1}\rho(t)\;dt,\quad (x>a)$$
and 
$$(I^{\alpha}_{b^-}\rho)(x):=\dis \frac{1}{\Gamma(\alpha)}\int_{x}^{b}(t-x)^{\alpha-1}\rho(t)\;dt,\quad (x<b),$$  
where $\Gamma$ denotes the usual Euler-Gamma function.
\end{definition}

Notice that if $\rho\in L^1(a,b)$, then by \cite[Theorem 1.5]{Baz},  the functions $I^{\alpha}_{a^+}\rho$ and $I^{\alpha}_{b^-}\rho$ belong to $W^{1,1}(a,b)$.

\begin{definition}\label{defDRL}
Let $\rho\in L^1(a,b)$ be such that the functions
\begin{align*}
x\mapsto \int_{a}^{x}(x-t)^{-\alpha}\rho(t)\;dt\;\mbox{ and } x\mapsto \int_{x}^{b}(t-x)^{-\alpha}\rho(t)\;dt
\end{align*}
belong to $W^{1,1}(a,b)$. The left,  and right Riemann–Liouville fractional derivatives of order $\alpha\in(0,1)$ of $\rho$,  are defined, respectively, by:
$$(\Dr^{\alpha}_{a^+}\rho)(x):=\frac{d}{dx}(I^{1-\alpha}_{a^+}\rho)(x)=\dis \frac{1}{\Gamma(1-\alpha)}\frac{d}{dx}
\int_{a}^{x}(x-t)^{-\alpha}\rho(t)\;dt,\quad (x>a)
$$
and
\begin{equation}\label{4}
(\D^{\alpha}_{b^-}\rho)(x):=-\frac{d}{dx}(I^{1-\alpha}_{b^-}\rho)(x)=\dis \frac{-1}{\Gamma(1-\alpha)}\frac{d}{dx}
\int_{x}^{b}(t-x)^{-\alpha}\rho(t)\;dt,\quad ~~(x<b).
\end{equation}
\end{definition}	

If $\alpha=1$ and $\rho\in W^{1,1}(a,b)$, then $(\Dr^{1}_{a^+}\rho)(x)=\rho^\prime(x)$ and $(\D^{1}_{b^-}\rho)(x)=-\rho^\prime(x)$.

\begin{definition}
The left, and right-sided Caputo fractional derivatives of order $\alpha\in(0,1)$ of $\rho\in W^{1,1}(a,b)$,  are defined respectively, by:
	\begin{equation}\label{4b1}
(\Dc^{\alpha}_{a^+}\rho)(x):=(I^{1-\alpha}_{a^+}\rho^\prime)(x)=\dis \frac{1}{\Gamma(1-\alpha)}\int_{a}^{x}(x-t)^{-\alpha}\rho^\prime(t)\;dt,\quad (x>a)
\end{equation}
and
\begin{equation}\label{4b}
(\Dc^{\alpha}_{b^-}\rho)(x):=-(I^{1-\alpha}_{b^-}\rho^\prime)(x)=\frac{-1}{\Gamma(1-\alpha)}\int_{x}^{b}(t-x)^{-\alpha}\rho^\prime(t)\;dt\quad (x<b).
\end{equation}
\end{definition}

As above,  if $\alpha=1$, then $(\Dc^{1}_{a^+}\rho)(x)=\rho^\prime(x)$ and $(\Dc^{1}_{b^-}\rho)(x)=-\rho^\prime(x)$.

We refer to the monograph \cite{Ga-Wa-2021} for the precise conditions on $\rho$ for which the integrals in \eqref{4b1}-\eqref{4b} exist.  We notice that the Caputo fractional derivative is meaningful for smooth functions $\rho$,  for any $0<\alpha\leq 1$, its existence requires the function $\rho$ to be absolutely continuous on $[a,b]$ which is equivalent to $\rho\in W^{1,1}(a,b)$.

Next, we shall just state some results that will be used in the present paper.  We start with the following lemma.   We refer to \cite[Theorem 3.5]{marichev} (see also \cite{kil}) for the proof.

\begin{lemma}\label{l1}
 Let $0<\alpha<1$,  $1<p<{1}/{\alpha}$, $q={p}/{(1-\alpha p)}$,  and $\rho\in L^p (a, b)$. Then, there is a constant $C=C(\alpha,p,q,a,b)>0$ such that 
 \begin{align*}	
\norm{I^{\alpha}_{a^+}\rho}_{L^q(a,b)}\leq&C\norm{\rho}_{L^p(a,b)}\\
\norm{I^{\alpha}_{b^-}\rho}_{L^q(a,b)}\leq&C\norm{\rho}_{L^p(a,b)}.
\end{align*}
\end{lemma}

\begin{remark}
Since the continuous embedding $L^2(a,b)\hookrightarrow L^p(a,b)$ holds for every $1\leq p\leq 2$, it follows from Lemma \ref{l1} that for every $0<\alpha<  1$ there is a constant $C>0$ such that for every $\rho\in L^2(a,b)$, 
\begin{equation}\label{marichev13bis}
\|I^\alpha_{a^+}\rho\|_{L^2(a,b)}\leq C\|\rho\|_{L^2(a,b)}.
\end{equation}
\end{remark}

Next, let $c_0,d_0\in \R$ and $0<\alpha\leq 1$. Let $\rho:[a,b]\to \R$ have the representation
\begin{align}\label{a1}	\rho(x)=\frac{c_0}{\Gamma(\alpha)}(x-a)^{\alpha-1}+I^{\alpha}_{a^+}\varphi(x)\;\mbox{ for a.e. } x\in [a,b],
	\end{align} and let also $g:[a,b]\to \R$ have the representation
\begin{align}\label{a2}
	g(x)=\frac{d_0}{\Gamma(\alpha)}(b-x)^{\alpha-1}+I^{\alpha}_{b^-}\psi(x)~~\mbox{ for a.e.  } x\in [a,b],
\end{align}
where $\varphi$ and $\psi$ belong to $L^2(a,b)$. We shall denote by $AC^{\alpha,2}_{a^+}$ and $AC^{\alpha,2}_{b^-}$ the spaces of all functions $\rho$ and $g$ having the representations \eqref{a1} and \eqref{a2}, respectively, with $\varphi , \psi \in L^2 (a, b)$.

\begin{remark}
Let $0<\alpha\leq 1$. We observe the following:
\begin{align}
	\D^{\alpha}_{a^+}\rho \in L^2 (a, b)&\iff \rho \in AC^{\alpha,2}_{a^+} \label{a3},\\ 
		\D^{\alpha}_{b^-}g\in L^2 (a, b)&\iff g \in AC^{\alpha,2}_{b^-} \label{a4}.
	\end{align}
\end{remark}

For more details  on these spaces and the proof of \eqref{a3}-\eqref{a4},  we refer to \cite{darius}.

We set
\begin{align}
 \Ha&=&AC^{\alpha,2}_{a^+} \cap L^2 (a, b)\label{defHa}\\
		\Hb&=&AC^{\alpha,2}_{b^-} \cap L^2 (a, b)\label{H2}.
	\end{align}
It follows from the definitions of $AC^{\alpha,2}_{a^+}$ and $AC^{\alpha,2}_{b^-}$ that,
\begin{align}
\rho \in \Ha \Longleftrightarrow \rho \in L^2(a, b) \hbox{ and } \D^{\alpha}_{a^+}\rho \in L^2(a, b),\label{prelim1}\\
	g \in \Hb \Longleftrightarrow g \in L^2(a, b) \text{ and } \D^{\alpha}_{b^-}g \in L^2(a, b).\label{prelim2}
\end{align}

We endow $\Ha$ with  the inner product
\begin{equation}\label{prodHa}
(\varphi,\psi)_{\Ha}=\int_a^b\varphi\psi \;dx+\int_a^b\D^{\alpha}_{a^+}\varphi\D^{\alpha}_{a^+}\psi\;dx.
\end{equation}
Then,  $\Ha$ endowed  with the norm given by
	\begin{equation}\label{normHa} \norm{\varphi}^2_{\Ha}=\norm{\varphi}^2_{L^2(a,b)}+\norm{\D^{\alpha}_{a^+}\varphi}^2_{L^2(a,b)}
	\end{equation}
is a Hilbert space  (see e.g. \cite{darius}).  Moreover,
the norm  on
	$\Ha$ given by \eqref{normHa} is equivalent to the norm defined by
	\begin{equation*} 
	\abs{\norm{\varphi}}^2=\abs{(I^{1-\alpha}
\varphi)(a)}^2+\norm{\D^{\alpha}_{a^+}\varphi}^2_{L^2(a,b)}~.
		\end{equation*}
In other words, there are two constants $0<C_1\leq C_2$  such that
		\begin{equation}\label{eq}
	C_1\abs{\norm{\varphi}}\leq	\norm{\varphi}_{\Ha}\leq C_2\abs{\norm{\varphi}}~~~\forall \varphi \in \Ha.
	\end{equation}
We refer to \cite[Theorem 2.4]{darius} for the proof of \eqref{eq}.

We also mention the following result that is interesting in its own but will not be used in the present paper. 
	
 \begin{lemma}\cite[Corollary 32]{darius} 
 Let $1/2<\alpha\leq 1$ . Then, the  continuous embedding
	\begin{equation} \label{a6}
	{\Ha} \hookrightarrow L^2(a,b)
	\end{equation}
	 is compact.
\end{lemma}

We mention that,  since throughout the paper we are not using the compact embedding \eqref{a6}, all our results hold true without any restriction on $\alpha$, that is, all the results hold for any $0<\alpha\leq 1$.
 
Next, let $0<\alpha\leq 1$ and set 
\begin{equation}\label{defVr}
\mathcal{V} := \Big\{y \in \Ha :\;  \Dc^{\alpha}_{b^-}(\beta \Dr^{\alpha}_{a^+}y) \in H_{b^-}^{1-\alpha}(a, b)\Big\},
\end{equation}
where $\beta\in C[a,b]$ and there is a constant $\beta_0>0$ such that  $\beta(x)\ge \beta_0$ for all $x\in[a,b]$.

From now on,  if $\mathbb{X}$ is a Banach space,  we shall denote by $L^2((0,T); \mathbb{X})$, the space of all measurable functions $f:(0,T)\to \mathbb{X}$ such that
$$\int_0^T \|f(t,\cdot)\|_{\mathbb{X}}^2\;dt<\infty.$$

We have the following result  that will be useful for some calculations in the upcoming sections.

\begin{lemma}\label{trace} 
Let $0<\alpha\leq 1$.
Then, the following assertions hold.
\begin{enumerate}[(i)]
\item Let $T>0$ and $\rho \in L^2((0,T);\Ha)$.  Then, for every $x_0\in [a,b],$ the function  $I_{a^+}^{1-\alpha}\rho(x_0,\cdot)$ exists   and belongs to $L^2(0,T)$. Moreover, there is a constant $C=C(a,b,\alpha)>0$ such that
  \begin{align}
    \|I_{a^+}^{1-\alpha}(\rho)(x_0,\cdot)\|^2_{L^2(0,T)}\leq C \|\rho\|^2_{L^2((0,T);\Ha)}. \label{traceI}
    \end{align}
  \item The space $\mathcal V$  defined in \eqref{defVr} endowed with the norm
\begin{align}\label{NM}
\|y\|_{\mathcal V}:=&\left(\|y\|_{\Ha}^2+\|  \Dc^{\alpha}_{b^-}(\beta \Dr^{\alpha}_{a^+}y)\|_{H_{b^-}^{1-\alpha}(a, b)}^2\right)^{1/2}\notag\\
=&\left(\|y\|_{\Ha}^2+ \|  I^{1-\alpha}_{b^-}(\beta \Dr^{\alpha}_{a^+}y)^\prime)\|_{L^2(a,b)}^2+\| (\beta\Dr_{a^+}^\alpha\rho)^\prime\|_{L^2(a, b)}^2\right)^{1/2},
\end{align}
and the associated scalar product
\begin{align}\label{SP}
(\phi,\psi)_{\mathcal V}:=&\int_a^b\phi\psi\;dx +\int_a^b \Dr_{a^+}^\alpha\phi \Dr_{a^+}^\alpha\psi\;dx +\int_a^b I^{1-\alpha}_{b^-}(\beta \Dr^{\alpha}_{a^+}\phi)^\prime I^{1-\alpha}_{b^-}(\beta \Dr^{\alpha}_{a^+}\psi)^\prime\;dx\notag\\
&+\int_a^b (\beta\Dr_{a^+}^\alpha \phi)^\prime (\beta\Dr_{a^+}^\alpha\psi)^\prime\;dx,
\end{align}
 is a Hilbert space.
 \item  Let $\rho \in \mathcal{V}$.  Then,  $\Dr_{a^+}^\alpha\rho\in C[a,b]$  and there is a constant $C=C(a,b,\alpha)>0$  such that 
 \begin{equation}\label{Trace} 
 \|\Dr_{a^+}^\alpha\rho\|_{C[a,b]}\leq C\|\rho\|_{\mathcal V}.
\end{equation}  
    \end{enumerate}
\end{lemma}

\begin{proof} 
The result is trivial if $\alpha=1$.  Now,  let $0<\alpha<1$.
We shall use the continuous embedding $H^1(a,b) \hookrightarrow C([a,b])$.

(i): Since $\rho \in L^2((0,T);\Ha)$,  we have that  $\rho\in L^2((0,T);L^2(a,b))$ and $\dis \Dr_{a^+}^\alpha\rho=\frac{d}{dx}\left(I_{a^+}^{1-\alpha}\rho\right) \in L^2((0,T); L^2(a,b))$. By Lemma \ref{l1}, $I_{a^+}^{1-\alpha}\rho \in L^2((0,T);L^2(a,b))$. Thus,  $I_{a^+}^{1-\alpha}\rho \in L^2((0,T);H^1(a,b))$.
Using the  continuous embedding $H^1(a,b)\hookrightarrow C([a,b])$,  we have that there is a constant $C=C(a,b)>0$ such that for  every $x_0\in [a,b]$,
          $$\begin{array}{lll}
          \dis \int_0^T \left|I_{a^+}^{1-\alpha}\rho(x_0,t)\right|^2\;dt&\leq &
          \dis \int_0^T \sup_{x\in [a,b]}\left|I_{a^+}^{1-\alpha}\rho(x,t)\right|^2\;dt\\
          &\leq &C
          \dis \int_0^T \|I_{a^+}^{1-\alpha}\rho (\cdot,t)\|^2_{H^1(a,b)}dt\\
          &=&\dis C\int_0^T \left(\|I_{a^+}^{1-\alpha}\rho(\cdot,t)\|^2_{L^2(a,b)}+
          \left\|\frac{d}{dx}I_{a^+}^{1-\alpha}\rho(\cdot,t)\right\|^2_{L^2(a,b)}\right)dt.
          \end{array}
          $$
          Observing that for a.e. $(x,t)\in (a,b)\times (0,T)$,
          \begin{equation*}
           \dis \frac{d}{dx}\left(I_{a^+}^{1-\alpha}\rho\right)(x,t)=\Dr_{a^+}^\alpha\rho(x,t),
           \end{equation*}
           we can deduce by using Lemma \ref{l1} again that
           \begin{align*}
          \dis \int_0^T \left|I_{a^+}^{1-\alpha}(\rho)(x_0,t)\right|^2dt\leq &
         \dis  C\int_0^T \left(\|I_{a^+}^{1-\alpha}\rho(\cdot,t)\|^2_{L^2(a,b)}+
          \left\|\frac{d}{dx}I_{a^+}^{1-\alpha}\rho(\cdot,t)\right\|^2_{L^2(a,b)}\right)\;dt\\
          \leq &
          C\left(\int_0^T \|\rho(\cdot,t)\|^2_{L^2(a,b)} \;dt+\int_0^T
          \left\|\Dr_{a^+}^\alpha\rho(\cdot,t)\right\|^2_{L^2(a,b)}\;dt\right)\\
          \leq &
          C\|\rho\|^2_{L^2((0,T);\Ha)} .
          \end{align*}
          This completes the proof of part (i).
          
   (ii): Since $\rho \in \mathcal{V}$,  we have that  $\Dc^\alpha_{b^-}(\beta\Dr_{a^+}^\alpha\rho)\in  H_{b^-}^{1-\alpha}(a,b)$. Thus,
       \begin{equation}\label{calcul1}
       \Dc^\alpha_{b^-}(\beta\Dr_{a^+}^\alpha\rho)\in L^2(a,b).
       \end{equation}
       and
       \begin{equation}\label{calcul2}
       \Dr_{b^-}^{1-\alpha}(\Dc^\alpha_{b^-}(\beta\Dr_{a^+}^\alpha\rho))\in L^2(a,b).
       \end{equation}
       Observe that on the one hand,  for a.e.  $x\in (a,b),$ we have
       $$
       \Dc^\alpha_{b^-}(\beta\Dr_{a^+}^\alpha\rho)(x)=
        I_{b^-}^{1-\alpha}\left(\frac{d}{dx}(\beta\Dr_{a^+}^\alpha\rho)(x)\right)= I_{b^-}^{1-\alpha}\left((\beta\Dr_{a^+}^\alpha\rho)^\prime(x)\right),$$
       and on the other hand, we have
       $$\dis
      \Dr_{b^-}^{1-\alpha} I_{b^-}^{1-\alpha}\left(\frac{d}{dx}(\beta\Dr_{a^+}^\alpha\rho)(x)\right)\\
      =
       \frac{d}{dx}(\beta\Dr_{a^+}^\alpha\rho)(x)=(\beta\Dr_{a^+}^\alpha\rho)^\prime(x).
       $$
       Notice that by definition
       \begin{equation}\label{ajout1}
       \Dc^\alpha_{b^-}(\beta\Dr_{a^+}^\alpha\rho)(x)= I_{b^-}^{1-\alpha}\left(\frac{d}{dx}(\beta\Dr_{a^+}^\alpha\rho)(x)\right)=I_{b^-}^{1-\alpha}\left((\beta\Dr_{a^+}^\alpha\rho)^\prime(x)\right).
       \end{equation}
       Since 
       \begin{equation}\label{ajout2}
      \Dr_{b^-}^{1-\alpha}(\Dc^\alpha_{b^-}(\beta\Dr_{a^+}^\alpha\rho))(x)=
       \frac{d}{dx}(\beta\Dr_{a^+}^\alpha\rho)(x)= (\beta\Dr_{a^+}^\alpha\rho)^\prime(x).
       \end{equation}
   Using the computations \eqref{ajout1} and \eqref{ajout2},  we can deduce that
\begin{align*}
\|  \Dc^{\alpha}_{b^-}(\beta \Dr^{\alpha}_{a^+}\rho)\|_{H_{b^-}^{1-\alpha}(a, b)}^2=& \|  \Dc^{\alpha}_{b^-}(\beta \Dr^{\alpha}_{a^+}\rho)\|_{L^2(a,b)}^2+ \|  \Dr_{b^-}^{1-\alpha}(\Dc^\alpha_{b^-}(\beta\Dr_{a^+}^\alpha\rho))\|_{L^2(a, b)}^2\\
=& \|  I^{1-\alpha}_{b^-}(\beta \Dr^{\alpha}_{a^+}\rho)^\prime)\|_{L^2(a,b)}^2+\| (\beta\Dr_{a^+}^\alpha\rho)^\prime\|_{L^2(a, b)}^2.
\end{align*}
It is clear that \eqref{NM} is a norm and \eqref{SP} is the associated scalar product. We show that $\mathcal V$ equipped with this norm is complete.  Indeed, let $\varphi_n\in\mathcal V$ be a Cauchy sequence with respect to the norm given in \eqref{NM}. Then,  $\varphi_n$ is a Cauchy sequence in $\Ha$. Since $\Ha$ is a Hilbert space,  it follows  that $\varphi_n$ converges to $\varphi\in \Ha$, as $n\to\infty$. We have to show that $(\beta\Dr_{a^+}^\alpha\varphi)^\prime, I^{1-\alpha}_{b^-}(\beta\Dr_{a^+}^\alpha\varphi)^\prime\in L^2(a,b)$,  and that $(\beta\Dr_{a^+}^\alpha\varphi_n)^\prime$ converges to $(\beta\Dr_{a^+}^\alpha\varphi)^\prime$ in $L^2(a,b)$,  $I^{1-\alpha}_{b^-}(\beta\Dr_{a^+}^\alpha\varphi_n)^\prime$ converges to $I^{1-\alpha}_{b^-}(\beta\Dr_{a^+}^\alpha\varphi)^\prime$ in $L^2(a,b)$, as $n\to\infty$. Since $(\beta\Dr_{a^+}^\alpha\varphi_n)^\prime$ is a Cauchy sequence in $L^2(a,b)$, we have that there is a function $\psi\in L^2(a,b)$ such that $(\beta\Dr_{a^+}^\alpha\varphi_n)^\prime\to \psi$ in $L^2(a,b)$,  as $n\to\infty$. Let $\phi\in C_c^\infty(a,b)$. Then,
\begin{align*}
\int_a^b(\beta\Dr_{a^+}^\alpha\varphi_n)^\prime\phi\;dx=-\int_a^b(\beta\Dr_{a^+}^\alpha\varphi_n)\phi^\prime\;dx.
\end{align*}
Noticing that $\beta\Dr_{a^+}^\alpha\varphi_n\to \beta\Dr_{a^+}^\alpha\varphi$ in $L^2(a,b)$ (since $\varphi_n\to\varphi\in \Ha$),  $\beta\in C[a,b]$,  and taking the limit of the latter identity as $n\to\infty$, we get that
\begin{align*}
\int_a^b\psi\phi\;dx=-\int_a^b(\beta\Dr_{a^+}^\alpha\varphi)\phi^\prime\;dx,\;\;\forall\;\phi\in C_c^\infty(a,b).
\end{align*}
It follows from the preceding identity that $\psi=(\beta\Dr_{a^+}^\alpha\varphi)^\prime\in L^2(a,b)$.   Now using Lemma \ref{l1}, we get that there is a constant $C>0$ such that
\begin{align*}
\|I^{1-\alpha}_{b^-}(\beta\Dr_{a^+}^\alpha\varphi_n)^\prime-I^{1-\alpha}_{b^-}(\beta\Dr_{a^+}^\alpha\varphi)^\prime\|_{L^2(a,b)}\leq C\|(\beta\Dr_{a^+}^\alpha\varphi_n)^\prime -(\beta\Dr_{a^+}^\alpha\varphi)^\prime\|_{L^2(a,b)}\to 0\;\mbox{ as } n\to\infty.
\end{align*}
We have shown that $\varphi\in\mathcal V$.

 (iii): Let $\rho \in \mathcal{V}$.  It follows from \eqref{calcul2} that
      \begin{equation}\label{mjw}
      (\beta\Dr_{a^+}^\alpha\rho)^\prime\in  L^2(a,b).
      \end{equation}
      Since $\rho \in H_{a^+}^{\alpha}(a,b)$ and $\beta\in C([a,b])$, it follows that  $\beta\Dr_{a^+}^\alpha\rho\in L^2(a,b)$.  Consequently,  using \eqref{mjw} we have that $\beta\Dr_{a^+}^\alpha\rho\in H^1(a,b)$.
      Since the embedding $H^1(a,b)\hookrightarrow C([a,b])$ is continuous, we have that   $(\beta\Dr_{a^+}^\alpha\rho)\in C[a,b]$.  Since $\beta\in C[a,b]$ and $\beta(x)\ge \beta_0>0$ for all $x\in [a,b]$, we can deduce that $\Dr_{a^+}^\alpha\rho\in C[a,b]$.  Thus, using \eqref{ajout2}, we can deduce that there is a constant $C>0$ such that
      \begin{align*}
      \|\Dr_{a^+}^\alpha\rho\|_{C[a,b]}\leq C\|\beta \Dr_{a^+}^\alpha\rho\|_{H^1(a,b)}\leq C\|\rho\|_{\mathcal V}.
      \end{align*}
 We have shown \eqref{Trace} and the proof is finished.
\end{proof}

\begin{remark}\label{remtrace}
Let $0<\alpha\leq 1$.
Notice that if $\rho\in L^2((0,T);\Ha)$,  then  from Lemma \ref{trace} we have that the traces
$(I^{1-\alpha}_{a^+})\rho)(a,\cdot)$ and $(I^{1-\alpha}_{a^+}\rho)(b,\cdot)$
 exist for a.e. $t\in (0,T)$, and belong to $L^2(0,T)$. Also if $\rho\in \mathcal V$, then $(\beta\Dr_{a^+}^\alpha\rho)(a)$ and $(\beta\Dr_{a^+}^\alpha\rho)(b)$ exists and are finite.
\end{remark}

Next, we introduce the following integration by parts formulas. We refer to \cite{Agr2007} for the proof.

\begin{lemma}\label{lem0}
Let $0<\alpha\leq 1$ and $y,\phi\in\Ha$. Then, the following assertions hold:
 \begin{equation}\label{integ1}
 \int_a^b\phi(s) \Dc_{b^-}^\alpha y(s)\;ds
 =-
 \left[y(s)(I_{a^+}^{1-\alpha}\phi)(s)\right]_{s=a}^{s=b}+\dis
 \int_a^b y(s) (\Dr_{a^+}^\alpha \phi)(s) \;ds
 \end{equation}
 and
 \begin{equation}\label{integ2}
 \int_a^b y(s) (\Dr_{a^+}^\alpha \phi)(s)y(s)\;ds
 =\left[y(s)(I_{a^+}^{1-\alpha}\phi)(s)\right]_{s=a}^{s=b}+\dis
 \int_a^b \phi(s) (\Dc_{b^-}^\alpha y)(s) \;ds.
 \end{equation}
\end{lemma}

More generally, we have the following.

\begin{lemma}\label{lem01}
Let $\beta\in C([a,b])$.  The the following assertions hold.
\begin{enumerate}
\item Le  $y,\phi\in \Ha$ be such that  $ \Dc_{b^-}^\alpha (\beta\Dr_{a^+}^\alpha y )\in L^2(a,b)$, and $(\beta\Dr_{a^+}^\alpha y )(b)$, $(\beta\Dr_{a^+}^\alpha y )(a)$ exist.
Then
 \begin{align}\label{integ1&2}
 \int_a^b\phi(s) \Dc_{b^-}^\alpha (\beta\Dr_{a^+}^\alpha y )(s)\;ds
 =-
 \left[(\beta\Dr_{a^+}^\alpha y )(s)I_{a^+}^{1-\alpha}(\phi)(s)\right]_{s=a}^{s=b}
 + \int_a^b (\beta\Dr_{a^+}^\alpha y )(s) \Dr_{a^+}^\alpha \phi(s)\;ds.
 \end{align}
 \item  Let $y,\phi\in \Ha$ be such that  $ \Dc_{b^-}^\alpha (\beta\Dr_{a^+}^\alpha y )$ and $ \Dc_{b^-}^\alpha (\beta\Dr_{a^+}^\alpha \phi )$ belong to $L^2(a,b)$, and $(\beta\Dr_{a^+}^\alpha y )(b)$, $(\beta\Dr_{a^+}^\alpha y )(a)$, $(\beta\Dr_{a^+}^\alpha \phi )(b)$, $(\beta\Dr_{a^+}^\alpha \phi )(a)$ exist. Then
 \begin{align}\label{IN-123}
 \int_a^b \Dc_{b^-}^\alpha\phi(s) (\beta\Dr_{a^+}^\alpha y )(s)\;ds =&-
 \left[(\beta\Dr_{a^+}^\alpha y )(s)I_{a^+}^{1-\alpha}(\phi)(s)\right]_{s=a}^{s=b}
 + \dis \left[I_{a^+}^{1-\alpha}(y)(s) \beta(s)\Dr_{a^+}^\alpha (\phi)(s)\right]_{s=a}^{s=b}\notag\\
 &+\dis
 \int_a^b y(s) \Dc_{b^-}^\alpha (\beta\Dr_{a^+}^\alpha (\phi))(s)\;ds.
 \end{align}
 \end{enumerate}
 \end{lemma}
 
 \begin{proof} 
 This follows directly from  \eqref{integ1} and \eqref{integ2}. 
  \end{proof}

 Throughout the remainder of the paper,  we assume the following.

 \begin{assumption}\label{asump-beta-q}
 We assume that $\beta\in C[a,b]$,  $q\in L^\infty(a,b)$, and there are two constants $\beta_0$ and $q_0$ such that
 \begin{align*}
\beta(x)\ge \beta_0>0 \; \mbox{ for all } x\in [a,b],\;\;\mbox{ and }  q(x)\ge q_0>0 \;\mbox{ for a.e.  }\; x\in [a,b].
 \end{align*}
 \end{assumption}

We have the following obvious result.

\begin{lemma}\label{lembil}
Let $0<\alpha\leq 1$, and $\beta$, $q$ satisfy Assumption \ref{asump-beta-q}.  For any $y,z\in \Ha$, we define the bilinear form $a(\cdot,\cdot):\Ha\times\Ha\to\R$ by:
\begin{align}\label{bil}
a(y,z)=\int_a^b \beta(x)\mathbb{D}^\alpha_{a^+}y(x) \mathbb{D}^\alpha_{a^+}z(x)\;dx +\int_a^b q(x)y(x)z(x)\;dx.
\end{align} Then, $a(\cdot,\cdot)$ is continuous and coercive. That is,  for every $y,\phi\in \Ha$, we have
 \begin{align}
 |a(y,\phi)|\leq \left(\|q\|_{\infty}+\|\beta\|_{\infty}\right)\|y\|_{\Ha}\|\phi\|_{\Ha},\label{bilest1}\\
 a(\phi,\phi)\geq \min\left(\beta_0,q_0\right)\|\phi\|^2_{\Ha}.\label{bilest2}
 \end{align}
In addition,  the form $a(\cdot,\cdot)$ is closed in $L^2(a,b)$.
\end{lemma}

We make the following observation.

\begin{remark}\label{rem-SG}
Let $0<\alpha\leq 1$. We have the following situation.
\begin{enumerate}
\item Let 
\begin{align}\label{defV0}
V_0:=\{u\in \Ha:\; (I_{a^+}^{1-\alpha})u(a)=0\}. 
\end{align}
Then, $V_0$ is a closed subspace of $\Ha$. Consider the bilinear, symmetric, closed, and coercive form given in \eqref{bil} but with  domain $D(a)=V_0$.  Let $A$ be the selfadjoint operator on $L^2(a,b)$ associated with $(a,V_0)$ in the sense that
\begin{equation*}
\begin{cases}
D(A)=\{y\in V_0:\; \exists\; f\in L^2(a,b)\;\mbox{ such that } a(y,z)=(f,z)_{L^2(a,b)}\;\forall\;z\in V_0\}\\
Ay=f.
\end{cases}
\end{equation*}
Using an integration by parts argument and the theory of distributions, one can show that
\begin{equation*}
\begin{cases}
D(A)=\{y\in V_0:\; \Dc^{\alpha}_{b^-}(\beta \Dr^{\alpha}_{a^+}y)+qy\in L^2(a,b)\;\mbox{ and } (\beta \Dr^{\alpha}_{a^+}y)(b^-)=0\}\\
Ay=\Dc^{\alpha}_{b^-}(\beta \Dr^{\alpha}_{a^+}y)+qy.
\end{cases}
\end{equation*}
The operator $-A$ generates a strongly continuous and analytic semigroup $S=(S(t))_{t\ge 0}$ in $L^2(a,b)$.
\item The operator $A$ can be also viewed as a bounded operator $\mathcal A:V_0\to V_0^\star$ given by
$$\langle \mathcal Ay,z\rangle_{V_0^\star,V_0}:=a(y,z),\;\forall \; y,z\in V_0\mbox{ with } \mathcal Ay=\Dc^{\alpha}_{b^-}(\beta \Dr^{\alpha}_{a^+}y)+qy.$$
The operator $-\mathcal A$ also generates a strongly continuous and analytic semigroup in $V_0^\star$ which coincides with the semigroup $S$ in $L^2(a,b)$. Here, $V_0^\star$ denotes the dual space of $V_0$ with respect to the pivot space $L^2(a,b)$, so that we have the continuous embeddings $V_0\hookrightarrow L^2(a,b)\hookrightarrow V_0^\star$.
\item Throughout the following,  if there is no confusion, we shall designate by $A$  both operators and $S$ the associated semigroups.
\end{enumerate}
\end{remark}

For notation convenience, let $V=\Ha$ and let $V^\star$ be its dual with respect to the pivot space $L^2(a,b)$.  As in Remark \ref{rem-SG},  for every $\phi\in V$, we have that $\Dc^{\alpha}_{b^-}(\beta \Dr^{\alpha}_{a^+}\phi)+q\phi$ belongs to $V^\star$. Using this fact we can rewrite Lemma \ref{lem01} as follows.

\begin{lemma}\label{lem01-2}
Let $0<\alpha\leq 1$,  $\beta\in C[a,b]$,  and $y,\phi\in V$ be such that $(\beta\Dr_{a^+}^\alpha y )(b)$, $(\beta\Dr_{a^+}^\alpha y )(a)$, $(\beta\Dr_{a^+}^\alpha \phi )(b)$, $(\beta\Dr_{a^+}^\alpha \phi )(a)$ exist.  Then
 \begin{align}\label{integ1&2-3}
\langle \Dc_{b^-}^\alpha (\beta\Dr_{a^+}^\alpha y ),\phi\rangle_{V^\star,V}
 =&-
 \left[(\beta\Dr_{a^+}^\alpha y )(s)(I_{a^+}^{1-\alpha}\phi)(s)\right]_{s=a}^{s=b}
 +
 \int_a^b (\beta\Dr_{a^+}^\alpha y )(s) (\Dr_{a^+}^\alpha \phi))(s)\;ds \notag\\
 =&-
 \left[(\beta\Dr_{a^+}^\alpha y )(s)(I_{a^+}^{1-\alpha}\phi)(s)\right]_{s=a}^{s=b}
 + \dis \left[(I_{a^+}^{1-\alpha}y)(s) (\beta\Dr_{a^+}^\alpha \phi)(s)\right]_{s=a}^{s=b}\notag\\
 &+\langle\Dc_{b^-}^\alpha (\beta\Dr_{a^+}^\alpha \phi),y\rangle_{V^\star,V} .
 \end{align}
 \end{lemma}
 
To conclude this section, we recall the following result taken from \cite[page 37]{lions1961} that will be useful in proving some existence results in the next sections.

\begin{theorem}\label{Theolions61} 
		Let $\left(F, \norm{\cdot}_F\right)$ be a Hilbert space. Let $\Phi$ be a subspace of $F$ endowed with a pre-Hilbert scalar product $(((\cdot,\cdot)))$  and associated norm $\||\cdot\||$.  Moreover, let $E:F\times \Phi\to \mathbb{C}$ be a sesquilinear form.  Assume that the following hypotheses hold:
		\begin{enumerate}
			\item  The embedding $\Phi \hookrightarrow F$ is continuous, that is, there is a constant $C_1>0$ such that
			\begin{equation}\label{theolions1}
			\norm{\varphi}_{F}\leq C_1|||\varphi|||~~\forall\; \varphi~ in~ \Phi.
			\end{equation}
			
			\item For all $\varphi\in \Phi$, the mapping $u\mapsto E(u,\varphi)$ is continuous on $F$.
			
			\item There is a constant $C_2>0$  such that
			\begin{equation}\label{theolions2}
			\abs{E(\varphi,\varphi)}\geq C_2 |||\varphi|||^2~~~\text{for all}~~\varphi\in \Phi.
			\end{equation}
		\end{enumerate}
		If $\varphi\mapsto\L(\varphi)$ is a continuous linear functional on $\Phi$, then there exists a function $u\in F$ verifying
		$$
		E(u,\varphi)=\L(\varphi)~~ \text{for all} ~~\varphi\in \Phi.
		$$
\end{theorem}

\section{Boundary optimal control problems on a single edge\label{oneedge}}
In this section, we are  concerned with the optimal control problem \eqref{OP1}-\eqref{ST1}.
Let 
\begin{align}\label{eq-32}
{J}(v):=\frac 12\int_Q \abs{y(v)-y_d}^2\,dxdt+\frac{N}{2}\int_0^T|v|^2 \;dt.
\end{align}
Here,  $y_d\in L^2(Q),$  $N>0,$  $\mathcal{U}_{ad}$   is a closed and convex subset of $L^2(0,T)$,  and $y=y(v)$ solves  the fractional Sturm–Liouville parabolic equation \eqref{ST1}.

To study the optimization problem \eqref{OP1}-\eqref{ST1},  we need some existence and regularity results of homogeneous and non-homogeneous fractional Sturm–Liouville parabolic equations in one dimension that we introduce in the next sections.

\subsection{Homogeneous fractional Sturm–Liouville parabolic equations in a single edge}

In this section, we are concerned with existence and regularity results of homogeneous fractional Sturm–Liouville parabolic equations of the type
\begin{equation}\label{pa}
\left\{\begin{array}{rlllllllll}
 y_t+\dis  \Dc_{b^-}^\alpha\,(\beta\, \Dr_{a^+}^\alpha y)+q\,y&=&f&\hbox{ in }& Q,\\
\dis (I_{a^+}^{1-\alpha} y)(a,\cdot)&=& 0&\hbox{ in }& (0,T),\\
\dis (\beta\Dr_{a^+}^{\alpha}y)(b^-,\cdot)&=&0&\hbox{ in }& (0,T),\\
y(\cdot,0)&=&y^0&\hbox{ in } &(a,b),
\end{array}
\right.
\end{equation}
that is, the system \eqref{ST1} with $v=0$.

From now on,  if $\mathbb{X} \hookrightarrow L^2(a,b)$ is a Hilbert space, then $\mathbb{X}^\star$ shall denote the dual of  $\mathbb{X}$  with respect to the pivot space $L^2(a,b)$, and by $\langle\cdot,\cdot\rangle_{\mathbb X^\star,\mathbb X}$ their duality map.

If we set
\begin{equation}\label{defW0T}
W(0,T;\mathbb{X}):= \left\{\zeta \in L^2((0,T);\mathbb{X}): \zeta_{t} \in L^2\left((0,T);\mathbb{X}^\star\right)\right\},
\end{equation}
then $W(0,T;\mathbb{X})$ endowed with the  norm given by
\begin{equation}\label{normW0T}
\|\psi\|^2_{W(0,T;\mathbb{X})}=\|\psi\|^2_{L^2((0,T);\mathbb{X})}+\|\psi_t\|^2_{
L^2\left((0,T);\mathbb{X}^\star\right)},\,\forall \psi \in W(0,T;\mathbb{X}),
\end{equation}
is a Hilbert space.  Moreover,  by \cite[Theorem II.5.12]{mag} we have the following continuous embedding:
\begin{equation}\label{contWTA}
W(0,T;\mathbb{X})\hookrightarrow C([0,T];L^2(a,b)).
\end{equation}

Throughout the following, without any mention, $V_0$ denotes the Hilbert space defined in \eqref{defV0} and $V_0^\star$ is its dual with respect to the pivot space $L^2(a,b)$.

Next, we introduce our notion of weak solutions to   the problem \eqref{pa}.

\begin{definition}\label{weaksolution}
Let  $0<\alpha\leq 1$,  $f\in L^2((0,T);V_0^\star)$ and $y^0\in L^2(a,b)$.  Let $\beta$, $q$ satisfy Assumption \ref{asump-beta-q}.  A function $y\in L^2((0,T);V_0)$ is said to be a weak solution of \eqref{pa} if the equality
\begin{equation}\label{weaksol}
\int_{0}^{T}\langle f,\,\varphi\rangle_{V_0^\star,V_0} \,dt +\int_a^b y^0(x)\varphi(x,0)\; dx=-\int_Qy \varphi_t\, dx\,dt+\int_{0}^{T}a(y,\varphi)dt,
\end{equation}
holds for every $\varphi \in  H(Q):=\Big\{\varphi\in L^2((0,T);V_0)\cap H^1((0,T);L^2(a,b)):\; \varphi(\cdot,T)=0 ~~\text{in}~~ (a,b) \Big\}$, 
  where the bilinear form $a(\cdot,\cdot)$ is defined in \eqref{bil}.
\end{definition}

We have the following existence and uniqueness result.

\begin{theorem}\label{Thm1}
Let  $0<\alpha\leq 1$,  $f\in L^2((0,T);(V_0)^\star)$,  and $y^0\in L^2(a,b)$.  Let $\beta$, $q$ satisfy Assumption \ref{asump-beta-q}.  Then, there exists a unique weak solution $y \in L^2((0,T);V_0)\cap H^1((0,T);V_0^\star)$ to \eqref{pa} in the sense of Definition~\ref{weaksolution}.
In addition, the following estimates hold true: There is a constant $C=C(\beta_0,q_0)>0$ such that
\begin{equation}\label{estimation1}
\begin{array}{lll}
\dis\|y(\cdot,T)\|^2_{L^2(a,b)}&\leq &C \left(\|y^0\|^2_{L^2(a,b)}+\|f\|^2_{L^2((0,T);V_0^\star)}\right)\\
\|y\|^2_{L^2((0,T);V_0)} &\leq &C
\left(\|y^0\|^2_{L^2(a,b)}+\|f\|^2_{L^2((0,T);V_0^\star)}\right)\\
\|y_t\|^2_{L^2((0,T);V_0^\star)} &\leq &C
\left(\|y^0\|^2_{L^2(a,b)}+\|f\|^2_{L^2((0,T);V_0^\star)}\right).
\end{array}
\end{equation}
\end{theorem}

\begin{proof} 
We proceed in four steps.

\textbf{Step 1.} We prove the existence by using Theorem \ref{Theolions61}. The norm on $L^2((0,T);V_0)$ is given by
$$\|\rho\|^2_{L^2((0,T);V_0)}=\int_0^T \|\rho(\cdot,t)\|^2_{V_0}\;dt=\int_0^T \|\rho(\cdot,t)\|^2_{\Ha}\;dt.$$
We consider the norm $\||\cdot\||$ on the pre-Hilbert space $H(Q)$ given by
$$
\||\rho\||^2:=\|\rho\|^2_{L^2((0,T);V_0)}+\|\rho(\cdot,0)\|^2_{L^2(a,b)},\, \forall \rho \in H(Q).
$$
It is clear that $\|\rho\|_{L^2((0,T);V_0)}\leq \||\rho\||$ for any $\rho\in H(Q)$.
This shows that we have the continuous embedding $H(Q)\hookrightarrow L^2((0,T);V_0)$.

Now, let $\varphi \in H(Q)$ and consider the bilinear form $\mathcal{F}(\cdot,\cdot)$ defined on $L^2((0,T);V_0)\times H(Q)$ by:
\begin{equation}\label{defCalF}
\mathcal{F}(y,\varphi)=-\int_Q  \varphi_{t}\,y\, dxdt+\int_0^T a(y,\varphi)\;dt.
\end{equation}
Using \eqref{bilest1}, we get that
$$\begin{array}{lllll}
\dis |\mathcal{F}(y,\varphi)|&\leq& \|y\|_{L^2(Q)}\|\varphi_{t}\|_{L^2(Q)}+ (\|q\|_{L^\infty(a,b)}+\|\beta\|_{L^\infty(a,b)})\|y\|_{L^2((0,T);V_0)}\|\varphi\|_{L^2((0,T);V_0)}\\
&\leq& \|y\|_{L^2((0,T);V_0)}\left(\|\varphi_{t}\|_{L^2(Q)}+ \left(\|q\|_{L^\infty(a,b)}+\|\beta\|_{L^\infty(a,b)}\right)\|\varphi\|_{L^2((0,T);V_0)}\right).
\end{array}
$$
This means that  there is a constant $C=C(\varphi,\|q\|_{L^\infty(a,b)},\|\beta\|_{L^\infty(a,b)})>0$  such that
$$|\mathcal{F}(y,\varphi)|\leq C\|y\|_{L^2((0,T);V_0)}. $$
Consequently, for every fixed $\varphi\in H(Q),$
the functional  $y\mapsto \mathcal{F}(y,\varphi)$ is continuous on $L^2((0,T);V_0).$

Next, using \eqref{bilest2}, we have that there is a constant $C>0$ such that for every  $\varphi\in H(Q)$,
\begin{align*}
\mathcal{F}(\varphi,\varphi)=&-\displaystyle\int_{0}^{T}\int_a^b \varphi \varphi_{t}\, dxdt+\dis\int_0^T a(\varphi,\varphi)\;dt\\
=&\frac 12  \|\varphi(\cdot,0)\|^2_{L^2(a,b)}+\int_0^T a(\varphi,\varphi)\;dt\\
\geq & C\||\varphi\||^2.
\end{align*}

Finally, let us consider the linear functional $L(\cdot):H(Q)\to \R$ defined by
$$
L(\varphi):=\dis \int_a^b y^0(x) \varphi(x,0)\, dx+\int_{0}^{T}\langle f,\, \varphi\rangle_{V_0^\star,V_0}\, dt .
$$
Then, there is a constant $C>0$ such that for every $\varphi\in H(Q)$,
\begin{align*}
|L(\varphi)|&\leq\|\varphi(\cdot,0)\|_{L^2(a,b)}\|y^0\|_{L^2(a,b)}+\|f\|_{L^2((0,T);V_0^\star)}
\|\varphi\|_{L^2((0,T);V_0)}\\
&\leq C\left(\|y^0\|_{L^2(a,b)}+\|f\|_{L^2((0,T);(V_0)^\star)}\right)\||\varphi\||.
\end{align*}
Therefore,  $L(\cdot)$ is continuous on $H(Q)$. It follows from Theorem \ref{Theolions61} that there exists $y\in L^2((0,T);V_0)$ such that
\begin{equation}\label{formvar}
\mathcal{F}(y,\varphi)= L(\varphi)\quad \forall \varphi \in H(Q).
\end{equation}
We have shown that the system \eqref{pa} has a solution $y\in L^2((0,T);V_0)$ in the sense of Definition \ref{weaksolution}.

 \textbf{Step 2.} We show that $y\in H^1((0,T);V_0^\star)$. It suffices to show that $y_t\in L^2((0,T);V_0^\star)$. Notice that $ \Dc_{b^-}^\alpha\,(\beta\, \Dr_{a^+}^\alpha y)+qy\in V_0^\star$ for a.e. $t\in (0,T)$ by Remark \ref{rem-SG}. Since $f\in L^2((0,T);V_0^\star)$,  we have that $y_t(\cdot,t)\in V_0^\star$ for a.e. $t\in (0,T)$.  It follows from \cite[Chapter IV, Section 1]{lions1961} (notice that it is a general result that applies to operators given by bilinear forms that may even dependent on the time $t$) that $y\in L^2((0,T);V_0)$ is a weak solution of \eqref{pa} in the sense of Definition \ref{weaksolution}, that is,  \eqref{weaksol} holds, if and only if the equality
 \begin{align*}
\frac{d}{dt} \left(y(\cdot,t),\varphi\right)_{L^2(\Omega} +a(y(\cdot,t),\varphi)=\langle f(\cdot,t),\varphi\rangle_{V_0^\star,V_0}
\end{align*}
holds for every $\varphi\in V_0$,  and a.e. $t\in (0,T)$.  Since $y_t(\cdot,t)\in V_0^\star$ for a.e. $t\in (0,T)$, it follows from the preceding identity that
\begin{align}\label{sol-et-dual}
 \langle y_t(\cdot,t),\varphi\rangle_{V_0^\star,V_0} +a(y(\cdot,t),\varphi)=\langle f(\cdot,t),\varphi\rangle_{V_0^\star,V_0},
\end{align}
for every $\varphi\in V_0$,  and a.e. $t\in (0,T)$.  
Thus,  using the continuity of the bilinear form $a$ (see \eqref{bilest1}),  taking $\varphi:=\phi(\cdot,t)$ (with $\phi\in L^2((0,T); V_0)$) as a test function in \eqref{sol-et-dual},  and integrating over $(0,T)$, we get that there is a constant $C>0$ (depending only on the coefficients of the operator) such that  for every $\phi\in L^2((0,T); V_0)$,
\begin{align*}
\left|\int_0^T\langle y_t,\phi\rangle_{V_0^\star,V_0}\;dt \right|\leq&C \int_0^T\|y(\cdot,t)\|_{V_0}\|\phi(\cdot,t)\|_{V_0}\;dt+\int_0^T\|f(\cdot,t)\|_{V_0^\star}\|\phi(\cdot,t)\|_{V_0}\;dt\\
\leq&C\|y\|_{L^2((0,T);V_0)}\|\phi\|_{L^2((0,T);V_0)}+\|f\|_{L^2((0,T);V_0^\star)}\|\phi\|_{L^2((0,T);V_0)}.
\end{align*}
Thus, $y_t\in L^2((0,T);V_0^\star)$ and 
\begin{equation}\label{mjw-2}
\|y_t\|_{ L^2((0,T);V_0^\star)}\leq C( \|y\|_{L^2((0,T);V_0)}+\|f\|_{L^2((0,T);V_0^\star)}). 
\end{equation}

 \textbf{Step 3.} We prove uniqueness.
Assume that there exist $y_1$ and $y_2$ two weak solutions of \eqref{pa} with the same right hand side $f$ and initial datum $y^0$.  Set $z:=y_1-y_2$. Then,  $z$ satisfies
\begin{equation}\label{pazero}
\left\{\begin{array}{rlllllllll}
z_t+\dis  \Dc_{b^-}^\alpha\,(\beta\, \Dr_{a^+}^\alpha z)+q\,z&=&0&\hbox{ in }& Q,\\
\dis (I_{a^+}^{1-\alpha} z)(a^+,\cdot)&=& 0&\hbox{ in }& (0,T),\\
\dis (\beta\Dr_{a^+}^{\alpha}z)(b^-,\cdot)&=&0&\hbox{ in }& (0,T),\\
z(\cdot,0)&=&0&\hbox{ in } &(a,b).
\end{array}
\right.
\end{equation}
Using  \eqref{sol-et-dual}  with $y=z$, taking $\varphi=z(\cdot,t)$ as a test function, and integrating over $(0,T)$, we get that
\begin{equation*}
0=\displaystyle\int_0^T\langle z_t, z\rangle_{V_0^\star,V_0}\,dt+\displaystyle\int_{0}^{T}a(z,z)\;dt.
\end{equation*}
Thus,  a simple integration gives
\begin{align*}
0=&\dis \frac 12\|z(\cdot,T)\|^2_{L^2(a,b)}+\int_0^T a(z,z)\;dt\\
=&\dis  \frac 12\|z(\cdot,T)\|^2_{L^2(a,b)}+\int_Q \beta(x)\left|\mathbb{D}^\alpha_{a^+}z(x,t)\right|^2\;dx dt+\int_Q q(x)|z(x,t)|^2\,dxdt\\
\geq&q_0\|z\|_{L^2(Q)}.
 \end{align*}
We can deduce that $z=0 $ in $Q$. Thus,  $y_1=y_2$ in $Q$ and we have shown uniqueness.

 \textbf{Step 4.} Finally, we show the estimates in \eqref{estimation1}.
Taking $\varphi=y(\cdot,t)$ as a test function  in \eqref{sol-et-dual}  and integrating over $(0,T)$, we obtain
\begin{align}\label{energie}
\dis  \frac 12\|y(\cdot,T)\|^2_{L^2(a,b)}+&\int_Q\beta(x)\left|\mathbb{D}^\alpha_{a^+}y(x,t)\right|^2\;dx dt+\int_Q q(x)|y(x,t)|^2\,dxdt\notag\\
=& \dis \frac 12 \|y^0\|^2_{L^2(a,b)}+\int_{0}^{T}\langle f,\, y\rangle_{(V_0)^\star,V_0}\, dt .
\end{align}
It follows from \eqref{bilest2} that
\begin{align*}
\int_a^b\beta(x)\left|\mathbb{D}^\alpha_{a^+}y(x)\right|^2 dx+\int_a^b q(x)|y(x)|^2\,dx\, \geq  \min(\beta_0,q_0)\|y\|^2_{V_0}.
\end{align*}
Thus, using Young's inequality with $\delta:=\min(\beta_0,q_0)>0$, we can deduce from \eqref{energie} that
\begin{align*}
  \frac 12\|y(\cdot,T)\|^2_{L^2(a,b)}+\delta\int_0^T\|y(\cdot,t)\|^2_{V_0}\;dt
\leq  \dis \frac{1}{2}\|y^0\|^2_{L^2(a,b)}+\frac{1}{2\delta}\|f\|^2_{L^2((0,T);V_0^\star)}+\frac{\delta}{2}
\|y\|^2_{L^2((0,T);V_0)}.
\end{align*}
The latter inequality implies that
\begin{align*}
 \frac 12 \|y(\cdot,T)\|^2_{L^2(a,b)}+\dis \frac{\delta}{2}\|y\|^2_{L^2((0,T);V_0)}
\leq  \dis \frac{1}{2}\|y^0\|^2_{L^2(a,b)}+\frac{1}{2\delta}\|f\|^2_{L^2((0,T);V_0^\star)}.
\end{align*}
Consequently,
$$\begin{array}{lll}
\dis\|y(\cdot,T)\|^2_{L^2(a,b)}&\leq &\dis \left(\frac{1}{\min(\beta_0,q_0)}+1\right) \left(\|y^0\|^2_{L^2(a,b)}+\|f\|^2_{L^2((0,T);V_0^\star)}\right),\\
\|y\|^2_{L^2((0,T);V_0)}
&\leq & \dis \left(\frac{1}{\min(\beta_0,q_0)}+\frac{1}{\left((\min(\beta_0,q_0)\right)^2}\right)
\left(\|y^0\|^2_{L^2(a,b)}+\|f\|^2_{L^2(0,T);V_0^\star}\right).
\end{array}
$$
Using the above two estimates and \eqref{mjw-2} we get that there is a constant $C>0$ such that
\begin{align*}
\|y_t\|^2_{L^2((0,T);V_0^\star)}
\leq  C\left(\|y^0\|^2_{L^2(a,b)}+\|f\|^2_{L^2(0,T);V_0^\star}\right).
\end{align*}
We have shown \eqref{estimation1} and the proof is finished.
\end{proof}

We conclude this section with the following important remark.

\begin{remark}\label{existy0T}
We make the following observations.
\begin{enumerate}
\item In Definition \ref{weaksolution},  we may replace $H(Q)$ with the space $\widetilde H(Q):= \{\varphi\in L^2((0,T);V_0)\cap H^1((0,T);V_0^\star):\;\varphi(\cdot,T)=0\mbox{ in } (a,b)\}$ and all the proofs remain the same (with the obvious modifications), but in that case we must replace \eqref{weaksol} with
\begin{equation}\label{weaksol-ww}
\int_{0}^{T}\langle f,\,\varphi\rangle_{V_0^\star,V_0} \,dt +\int_a^b y^0(x)\varphi(x,0)\; dx=-\int_0^T\langle  \varphi_t,y\rangle_{V_0^\star,V_0}\,dt+\int_{0}^{T}a(y,\varphi)dt,
\end{equation}
holds for every $\varphi\in \widetilde H(Q)$.
\item Notice that if a weak solution $y$ of \eqref{pa} exists, then it belongs to $W(0,T;V_0)$ (see \eqref{defW0T}).  Thus, it follows from \eqref{contWTA} that $y(\cdot,0)$ and $y(\cdot,T)$ exist and belong to $L^2(a,b)$. 
\item Using the operator $A$ defined in Remark \ref{rem-SG}, we have that the system \eqref{pa} can be rewritten as the following abstract Cauchy problem:
\begin{equation}\label{ACP}
\begin{cases}
y_t +Ay= f \;&\mbox{ in }\ Q\\
y(0)= y^0 &\mbox{ in }(a,b).
\end{cases}
\end{equation}
It follows from semigroups theory that for every $f\in L^2((0,T);V_0^\star)$ and $y_0\in L^2(a,b)$,  the Cauchy problem \eqref{ACP} (and hence, the system \eqref{pa}) has a unique weak solution $y$ given by
\begin{align*}
y(\cdot,t)=S(t)y^0+\int_0^t S(t-\tau)f(\cdot,\tau)\;d\tau
\end{align*}
where we recall that $S=(S(t))$ is the semigroup generated by the operator $-A$.  In the previous proof we have chosen to use Theorem \ref{Theolions61} since this method will be crucial for the non-homogeneous problem in the next sections where semigroups theory cannot be applied directly.
\item If $f\in L^2(Q)$ and $y^0\in D(A)$, then the regularity of $y$ can be improved. In fact, in that case, using maximal regularity results for abstract Cauchy problems (see e.g. \cite{de-Simon}), we can deduce that $y$ is a strong solution, that is, $y\in C([0,T];D(A))\cap H^1((0,T);L^2(a,b))$ and the first equation in  \eqref{pa} is satisfied pointwise for a.e.  $(x,t)\in Q$, where $D(A)$ is the domain of the operator $A$ defined in Remark \ref{rem-SG}(a).
\end{enumerate}
\end{remark}

\subsection{Non-homogeneous fractional Sturm–Liouville parabolic equations  in a single edge}

Now we are interested in the non-homogeneous problem \eqref{ST1}.  Recall the Hilbert space $V_0$ given in \eqref{defV0}.
Then, we have the following continuous embeddings:
 \begin{equation}\label{injection1a}
V_0\hookrightarrow H_{a^+}^{\alpha}(a,b)\hookrightarrow L^2(a,b)\hookrightarrow \left(H_{a^+}^{\alpha}(a,b)\right)^\star\hookrightarrow V_0^\star.
\end{equation}
Throughout this section, $H(Q)$ denotes the space introduced in Definition \ref{weaksolution}.

Here,  is our notion of solutions.

 \begin{definition}\label{weaksolutionnh}
Let  $0<\alpha\leq 1$,  $f\in L^2((0,T);V_0^\star),$  $y^0\in L^2(a,b)$,  and $v\in L^2(0,T)$. Let also $q$ and $\beta$ satisfy Assumption \ref{asump-beta-q}.  A function $
y\in L^2((0,T);V_0)$ is said to be a weak solution of  \eqref{ST1} if the  equality
\begin{equation}\label{weaksolnh}
\mathcal{F}(y,\varphi)=L_0(\phi), \quad \forall \phi\in H(Q)
\end{equation}
holds,  where $\mathcal{F}(\cdot,\cdot)$ is given by  \eqref{defCalF} and  $L_0(\cdot)$ is defined by
$$
L_0(\phi)=\dis \int_a^b y^0(x)\phi(x,0)\, dx+\int_{0}^{T}\langle f,\phi\rangle_{V_0^\prime,V_0}\, dt+\int_0^T v(t)(I^{1-\alpha}_{a^+}\phi) (b^-,t)\;dt .
$$
\end{definition}

Recall that the pre-Hilbert space $H(Q)$ is endowed with the norm given by
$$
\||\rho\||^2=\|\rho\|^2_{L^2((0,T);V_0)}+\|\rho(\cdot,0)\|^2_{L^2((a,b))},\, \forall \rho \in H(Q).
$$

We have the following existence and uniqueness result.

 \begin{theorem}\label{theoexistnh}
 Let  $0<\alpha\leq 1$,   $v\in L^2(0,T)$,  and $f\in L^2((0,T);V_0^\star)$.  Let $\beta$, $q$ satisfy Assumption \ref{asump-beta-q}.   Then, there exists a unique weak solution $y\in L^2((0,T);V_0)\cap H^1((0,T);V_0^\star)$ of \eqref{ST1} in the sense of Definition \ref{weaksolutionnh}. Moreover, the following estimates hold true: There is a constant $C=C(a,b,\beta_0,q_0)>0$ such that 
\begin{align}\label{estimation2}
 \|y(\cdot,T)\|^2_{L^2(a,b)}\leq & C\left( \|y^0\|^2_{L^2(a,b)}+\|f\|^2_{L^2((0,T);V_0^\star)}+\|v\|^2_{L^2(0,T)}\right)\\
\|y\|^2_{L^2((0,T);V_0)}
\leq & C\left(\|y^0\|^2_{L^2(a,b)}+\|f\|^2_{L^2((0,T);V_0^\star)}+\|v\|^2_{L^2(0,T)}\right)\notag\\
\|y_t\|^2_{L^2((0,T);V_0^\star)}
\leq & C\left(\|y^0\|^2_{L^2(a,b)}+\|f\|^2_{L^2((0,T);V_0^\star)}+\|v\|^2_{L^2(0,T)}\right).\notag
\end{align}
\end{theorem}

\begin{proof} We proceed as in the above homogeneous case by using Theorem \ref{Theolions61} again.

 \textbf{Step 1.} We prove existence. 
Recall that we have the continuous embedding $H(Q)\hookrightarrow L^2((0,T);V_0)$.
It follows from the computation  in the homogeneous case that for every $ \varphi\in H(Q)$, the functional $y\mapsto \mathcal{F}(y,\varphi)$ is continuous on $L^2((0,T);V_0)$.  More precisely, there is a constant $C>0$ such that for all $(y,\varphi)\in L^2((0,T);V_0)\times H(Q)$, we have
$$|\mathcal{F}(y,\varphi)|\leq C\|y\|_{L^2((0,T);V_0)}\||\varphi\||. $$
On the other hand,  $\mathcal{F}(\cdot,\cdot)$ is coercive on $H(Q)$, since for every $\varphi\in H(Q)$ we have
$$
\mathcal{F}(\varphi,\varphi)\geq
\min\left(\frac{1}{2},\min(\beta_0,q_0)\right)\||\varphi\||^2.
$$
To complete the prove of the existence, we need to show that the linear functional $L_0(\cdot)$ is continuous on $H(Q)$.
Using Cauchy-Schwarz's inequality, we can estimate
\begin{equation}\label{observation1}\begin{array}{llllll}
\left|(I^{1-\alpha}_{a^+}\varphi)(b^-,t)\right|&=& \left|(I^{1-\alpha}_{a^+}\varphi)(b^-,t)-(I^{1-\alpha}_{a^+}\varphi)(a,t)\right|=\dis \left|\int_a^b\frac{d}{dx}(I^{1-\alpha}_{a^+}\varphi)(x,t)\;dx\right| \\
&=&	\dis \left|\int_a^b\D^{\alpha}_{a^+}\varphi(x,t)\;dx\right| \leq	 \abs{b-a}^{1/2}\left\|\D^{\alpha}_{a^+}\varphi(\cdot,t)\right\|_{L^2(a,b)}.
\end{array}
\end{equation}
Thus, there is a positive constant $C=C(y^0,f,v,b,a)$ such that
\begin{align*}
|L_0(\varphi)|\leq&\|\varphi(\cdot,0)\|_{L^2(a,b)}\|y^0\|_{L^2(a,b)}+\|f\|_{L^2((0,T);V_0^\star)}
\|\varphi\|_{L^2((0,T);V_0)}+\dis \int_0^T |v(t)|\left| (I^{1-\alpha}_{a^+}\varphi) (b^-,t)\right|\;dt \\
\leq&\|\varphi(\cdot,0)\|_{L^2(a,b)}\|y^0\|_{L^2(a,b)}+\|f\|_{L^2((0,T);V_0^\star)}
\|\varphi\|_{L^2((0,T);V_0)}\\
&+\dis \abs{b-a}^{1/2}\int_0^T |v(t)|\left\|\D^{\alpha}_{a^+}\varphi(\cdot,t)\right\|_{L^2(a,b)}\;dt \\
\leq&C\left(\|\varphi(\cdot,0)\|^2_{L^2(a,b)}+\|\varphi\|^2_{L^2((0,T);V_0)}+
\left\|\D^{\alpha}_{a^+}\varphi\right\|^2_{L^2(Q)}\right)^{1/2}\\
\leq&C\||\varphi\||.
\end{align*}
We have shown that $L_0(\cdot)$ is continuous on $H(Q)$. It follows from Theorem \ref{Theolions61} again that there exists $y\in L^2((0,T);V_0)$ satisfying \eqref{weaksolnh}.   In addition, by Remark \ref{rem-SG}, we have that $y_t(\cdot,t)\in V_0^\star$ for a.e. $t\in (0,T)$.  As above,  it follows from \cite[Chapter IV, Section 1]{lions1961} that $y$ is a weak solution of \eqref{ST1} in the sense of Definition \ref{weaksolutionnh}, that is, \eqref{weaksolnh} holds, if and only if, the equality
\begin{align*}
\frac{d}{dt}(y(\cdot,t),\varphi)_{L^2(\Omega)}+a(y(\cdot,t),\varphi)=\langle f(\cdot,t), \varphi\rangle_{V_0^\star,V_0}
+v(t)(I_{a^+}^{1-\alpha}\varphi)(b^-)
\end{align*}
holds for every $\varphi\in V_0$, and a.e. $t\in (0,T)$.  Since $y_t(\cdot,t)\in V_0^\star$ for a.e. $t\in (0,T)$,  we can deduce from the preceding equality that 
\begin{align}\label{w-sol-nh}
\frac{d}{dt}\langle y_t(\cdot,t),\varphi\rangle_{V_0^\star,V_0}+a(y(\cdot,t),\varphi)=\langle f(\cdot,t), \varphi\rangle_{V_0^\star,V_0}
+v(t)(I_{a^+}^{1-\alpha}\varphi)(b^-)
\end{align}
holds for every $\varphi\in V_0$, and a.e. $t\in (0,T)$. 
Hence, proceeding as in Step 2 in the proof of Theorem \ref{Thm1} by taking $\varphi=\phi(\cdot,t)$ (with $\phi \in L^2((0,T);V_0)$) in \eqref{w-sol-nh} and integrating over $(0,T)$,  we get that if $y\in L^2((0,T);V_0)$ and satisfies \eqref{ST1} with $f\in L^2((0,T);V_0^\star)$,  then $ y_t\in L^2((0,T);V_0^\star)$. Thus,  $y\in H^1((0,T);V_0^\star)$. We have also shown that
$y\in W(0,T;V_0)$.  Thus, $y(\cdot,0)$ and $y(\cdot,T)$ exist and belong to $L^2(a,b)$. 

 \textbf{Step 2.} We show uniqueness.
Assume that \eqref{ST1} has two solutions $y_1$ and $y_2$  with the same right hand side data $f$, $v$, and the same initial datum $y^0$.  Set $z:=y_1-y_2$. Then,  $z$ satisfies \eqref{pazero}.  It follows from Step 3 in the proof of Theorem \ref{Thm1} 
 that $z=0 $ in $Q$. Thus,  $y_1=y_2$  in $Q$ and we have shown uniqueness.

 \textbf{Step 3.} Finally,  we show the estimates in \eqref{estimation2}.
Taking $\varphi=y(\cdot,t)$ as a test function in \eqref{w-sol-nh} and integrating over $(0,T)$,  we get
\begin{align}\label{energie1}
 \frac 12 \|y(\cdot,T)\|^2_{L^2(a,b)}+&\int_Q\beta(x)\left|\mathbb{D}^\alpha_{a^+}y(x,t)\right|^2 dx \, dt+\int_Q q(x)|y(x,t)|^2\,dx\, dt\notag\\
=&
  \frac 12\|y^0\|^2_{L^2(a,b)}+\int_{0}^{T}\langle f,\, y\rangle_{V_0^\star,V_0}\, dt +\int_0^T v(t)(I^{1-\alpha}_{a^+}y) (b^-,t)\;dt.
\end{align}
Since
\begin{align*}
\int_a^b\beta(x)\left|\mathbb{D}^\alpha_{a^+}y(x)\right|^2 dx+\int_a^b q(x)|y(x)|^2\,dx\, \geq \min(\beta_0,q_0)\|y\|^2_{V_0},
\end{align*}
using Young's inequality, \eqref{observation1},  and the latter estimate,  we can deduce from \eqref{energie1} that
\begin{align*}
\dis  \frac 12\|y(\cdot,T)\|^2_{L^2(a,b)}+\min(\beta_0,q_0)\int_0^T\|y\|^2_{V_0}dt
\leq & \dis  \frac 12\|y^0\|^2_{L^2(a,b)}+\frac{1}{2\delta}\|f\|^2_{L^2((0,T);(V_0)^\star)}+\frac{\delta}{2}
\|y\|^2_{L^2((0,T);V_0)}\\
&+\frac{b-a}{2\delta}\|v\|^2_{L^2(0,T)}+
\frac{\delta}{2}\left\|\mathbb{D}^\alpha_{a^+}y\right\|^2_{L^2(Q)}\\
\leq & \dis  \frac 12\|y^0\|^2_{L^2(a,b)}+\frac{1}{2\delta}\|f\|^2_{L^2((0,T);(V_0)^\star)}+\frac{\delta}{2}
\|y\|^2_{L^2((0,T);V_0)}\\
&+\frac{b-a}{2\delta}\|v\|^2_{L^2(0,T)}+
\frac{\delta}{2}\left\|y\right\|^2_{L^2((0,T);V_0)},
\end{align*}
for every $\delta>0$.  Choosing  $\delta:= \frac{\min(\beta_0,q_0)}{2}>0$, we can deduce from the latter inequality that
\begin{align*}
\dis \|y(\cdot,T)\|^2_{L^2(a,b)}&\leq  \left(1+ \frac{2(b-a+1)}{\min(\beta_0,q_0)} \right)\left( \|y^0\|^2_{L^2(a,b)}+\|f\|^2_{L^2((0,T);V_0^\star)}+\|v\|^2_{L^2(0,T)}\right),\\
\|y\|^2_{L^2((0,T);V_0)}
&\leq  \dis \left(\frac{1}{\min(\beta_0,q_0)} + \frac{2(b-a+1)}{\min(\beta_0,q_0)^2}\right)\left(\|y^0\|^2_{L^2(a,b)}+\|f\|^2_{L^2((0,T);V_0^\star)}+\|v\|^2_{L^2(0,T)}\right).
\end{align*}
We have shown the first two estimates in \eqref{estimation2}.  The last one can be obtained by proceeding exactly as in the proof of \eqref{mjw-2}.  The proof is finished.
\end{proof}

\begin{remark}\label{remf}
Notice that Theorems \ref{Thm1} and \ref{theoexistnh} trivially hold if $f\in L^2(Q)\hookrightarrow L^2((0,T);V_0^\star)$.
\end{remark}

\subsection{Existence of minimizers and optimality conditions in a single edge}

In this section we  are concerned with the existence of minimizers of the optimal control problem \eqref{OP1}-\eqref{ST1}.  Recall that $\mathcal{U}_{ad}$ is a closed and convex subset of $L^2(0,T)$. We have the following result.

\begin{theorem}\label{theoopt}
Let  $0<\alpha\leq 1$.   Let $\beta$, $q$ satisfy Assumption \ref{asump-beta-q}.  Then, there exists a unique solution $u\in \mathcal{U}_{ad}$ of the optimal control problem \eqref{OP1}-\eqref{ST1}.
\end{theorem}

\begin{proof} Let $y(0)$ satisfy \eqref{ST1} with $v=0.$ Then,  in view of  \eqref{ST1}, $z=z(v):=y(v)-y(0)$ satisfies
$$\left\{\begin{array}{rlllllllll}
z_t+\dis  \Dc_{b^-}^\alpha\,(\beta\, \Dr_{a^+}^\alpha z)+q\,z&=&0&\hbox{ in }& Q,\\
\dis (I_{a^+}^{1-\alpha} z)(a,\cdot)&=& 0&\hbox{ in }& (0,T),\\
\dis (\beta\Dr_{a^+}^{\alpha}z)(b^-,\cdot)&=&v&\hbox{ in }& (0,T),\\
z(\cdot,0)&=&0&\hbox{ in } &(a,b).
\end{array}
\right.
$$
Since $v\in L^2(0,T)$, it follows from \eqref{estimation2} that
\begin{equation}\label{est0}
\dis\|z\|_{L^2((0,T);V_0)}
\leq  \dis \sqrt{\left(\frac{1}{\min(\beta_0,q_0)} + \frac{2(b-a+1)}{\min(\beta_0,q_0)^2}\right)}\|v\|_{L^2(0,T)}.
\end{equation}
Now, let $\pi:\mathcal{U}_{ad}\times \mathcal{U}_{ad}\to\mathbb R$ be given by
$$\begin{array}{llll}
 \pi(u,v)&:=&\dis  \int_Q (y(v)-y(0))(y(u)-y(0))\, dx dt+ N\int_0^T  uv\; dt\\
 &=&\dis  \int_Q z(v)z(u)\, dxdt+ N\int_0^T  uv\; dt.
 \end{array}$$
Then,
$\pi(\cdot,\cdot)$ is continuous, because for every $u,v\in \mathcal{U}_{ad}$ we have
\begin{align*}
|\pi(u,v)|&=\dis \left |\int_Q (y(v)-y(0))(y(u)-y(0))\, dx dt+ N\int_0^T  uv\; dt\right|\\
&=\dis \left |\int_Q z(v)z(u)\, dx \, dt+ N\int_0^T  uv \;dt\right|\\
&\leq\|z(v)\|_{L^2(Q)}\|z(u)\|_{L^2(Q)}+N\|v\|_{L^2(0,T)}\|u\|_{L^2(0,T)}\\
&\leq\left(\frac{1}{\min(\beta_0,q_0)} + \frac{2(b-a+1)}{\min(\beta_0,q_0)^2}\right)\|v\|_{L^2(0,T)}\|u\|_{L^2(0,T)}+N\|v\|_{L^2(0,T)}\|u\|_{L^2(0,T)}\\
&\leq\left[\left(\frac{1}{\min(\beta_0,q_0)} + \frac{2(b-a+1)}{\min(\beta_0,q_0)^2}\right)+N\right]\|v\|_{L^2(0,T)}\|u\|_{L^2(0,T)}.
\end{align*}
In addition,  $\pi(\cdot,\cdot)$ is coercive, because for every $v\in\mathcal{U}_{ad}$ we have
$$\pi(v,v)=\dis \int_Q |y(v)-y(0)|^2\, dx \, dt+ N\int_0^T  |v|^2\;dt\geq N \|v\|^2_{L^2(0,T)}.$$
Observe that the functional $$v\mapsto \dis \int_Q (y(v)-y(0))(y_d-y(0)) \, dxdt$$
 is  linear and continuous on $\mathcal{U}_{ad}$, because
$$\begin{array}{lllll}
\dis \left|\int_Q (y(v)-y(0))(y_d-y(0)) \, dx\, dt\right|&=&\dis
\left|\int_Q z(v)(y_d-y(0)) \, dx\, dt\right|\\
&\leq &\dis \|y_d-y(0)\|_{L^2(Q)}\|z(v)\|_{L^2(Q)}\\
&\leq & \dis \|y_d-y(0)\|_{L^2(Q)}\sqrt{\left(\frac{1}{\min(\beta_0,q_0)} + \frac{2(b-a+1)}{\min(\beta_0,q_0)^2}\right)}\|v\|_{L^2(0,T)}.
\end{array}$$
On the other hand,  since the cost functional can be rewritten as
$$J(v)= \frac 12\left[\pi(v,v)-2\int_Q (y(v)-y(0))(y_d-y(v)) \, dx \,dt+ \|y(0)-y_d\|^2_{L^2(Q)}\right],$$
using  \cite[Theorem 1.1 Page 4]{lions1968}, we can deduce  that there exists a unique $u\in \mathcal{U}_{ad}$ solution of the minimization problem \eqref{OP1}-\eqref{ST1}. The proof is concluded.
\end{proof}

The following result gives a characterization of the optimal control.

\begin{theorem}\label{theo-39}
Let  $0<\alpha\leq 1$.   Let $\beta$, $q$ satisfy Assumption \ref{asump-beta-q}.
	Let $u\in \mathcal{U}_{ad}$ be the optimal control for the problem \eqref{OP1}-\eqref{ST1}. Then,  the first order necessary optimality conditions are given by
	\begin{equation}\label{pa6}
	\int_0^T\Big( Nu(t) - (I^{1-\alpha}_{a^+}p)(b^-,t)\Big)(v(t)-u(t))\;dt\geq 0 ~~\forall v\in \mathcal{U}_{ad},
	\end{equation}	
	where   $p$ is the unique weak (strong) solution of the adjoint system
	\begin{equation}\label{pa5}
\left\{
\begin{array}{llllllllllll}
\displaystyle	- p_t+\dis  \Dc_{b^-}^\alpha\,(\beta\, \Dr_{a^+}^\alpha p)+q\,p&=&y_d-y&\text{in}& Q,\\
\displaystyle	(I^{1-\alpha}_{a^+}p)(a,\cdot)&=&0&\text{in}&  (0,T),\\
\displaystyle	(\beta\D^{\alpha}_{a^+}p)(b^-,\cdot)&=&0&\text{in}& (0,T),\\
\displaystyle	p(\cdot,T)&=&0 &\text{in}& (a,b),
	\end{array}
	\right.
\end{equation}
and $y$ solves the state equation \eqref{ST1} with $v=u$.   In addition, \eqref{pa6} is equivalent to
\begin{equation}\label{scpr}
u=\mathcal P_{\mathcal U_{ad}}\Big(N^{-1}(I_{a^+}^{1-\alpha}p)(b^-,\cdot)\Big),
\end{equation}
where $\mathcal P_{\mathcal U_{ad}}$ is the projection onto the set $\mathcal U_{ad}$. Moreover, since our functional $J$ is convex, we have that \eqref{pa6} is also a sufficient condition.
\end{theorem}

\begin{proof}
We write the Euler-Lagrange first order optimality conditions that characterize the optimal control $u$. That is,
\begin{equation}\label{euler1}
\dis \lim_{\theta \downarrow 0} \frac{J(u+\theta(v-u))-J(u)}{\theta} \geq 0,\quad \forall v\in \mathcal{U}_{ad}.
\end{equation}
After some straightforward  calculations, \eqref{euler1} gives
\begin{equation}\label{euler2}
\int_Q z(y-y_d)\, dx\, dt+N\int_0^T u\, (v-u)\;dt \geq 0 ,\quad \forall v\in \mathcal{U}_{ad},
\end{equation}
where $z=z(v-u)\in L^2((0,T);V_0)\cap H^1((0,T);V_0^\star)$ is the unique weak solution of
\begin{equation}\label{pa4z}
\left\{
\begin{array}{llllllllllll}
\displaystyle z_t+\dis  \Dc_{b^-}^\alpha\,(\beta\, \Dr_{a^+}^\alpha z)+q\,z&=&0 &\text{in}& Q,\\
\displaystyle (I^{1-\alpha}_{a^+}z)(a,\cdot)&=&0&\text{in}&  (0,T),\\
\displaystyle (\beta\D^{\alpha}_{a^+}z)(b^-,\cdot)&=&v-u&\text{in}&(0,T),\\
\displaystyle	z(\cdot,0)&=&0 &\text{in}& (a,b).
	\end{array}
	\right.
	\end{equation}
Notice that it follows from Remark \ref{existy0T} that $p\in L^2((0,T);V_0)\cap H^1((0,T);L^2(a,b))$, and $p$ is a strong solution.  Hence,  we can deduce that $\Dc_{b^-}^\alpha\,(\beta\, \Dr_{a^+}^\alpha z)+q\in L^2(Q)$.
Multiplying  the first equation in \eqref{pa4z} by $p$ strong solution of \eqref{pa5}, we obtain after integration by parts,
 \begin{align*}
 0=\int_Qz(-p_t+\dis  \Dc_{b^-}^\alpha\,(\beta\, \Dr_{a^+}^\alpha p)+q\,p)\;dx dt- \int_0^T (v-u)(I_{a^+}^{1-\alpha}p)(b^-,t)\; dt.
 \end{align*}
Thus,
\begin{equation}\label{mmm}
  0= \int_Q z(y_d-y)\;dxt dt- \int_0^T (v-u)(I_{a^+}^{1-\alpha}p)(b^-,t)\; dt.
  \end{equation}
Combining \eqref{euler2}-\eqref{mmm},  we get \eqref{pa6}.  Finally,  that \eqref{pa6} is equivalent to
\eqref{scpr} follows directly form \cite[Theorem 3.3.5]{Att-BM}.
The proof is concluded.
\end{proof}

	\section{Boundary optimal  control problems in general star graphs}\label{graph}

In this section,  we consider the fractional Sturm-Liouville problem in a general star graph,  that is, the system \eqref{ST2}.
Recall that $0<\alpha\leq 1$,  $f=(f^i)\in L^2((0,T);\mathbb{L}^2)$, the control functions  $u_i\in L^2(0,T),\, i=2,\ldots,  m$, and  $v_i\in L^2(0,T),\, i=m+1,\ldots, n,$  $Q_i:= (a,b_i)\times(0,T),\,i=1,\ldots, n$, where we have set $\dis\mathbb{L}^2:=\dis \prod_{i= 1}^n L^2(a,b_i)$.

Recalling the  formulation with non-homogeneous Dirichlet
boundary conditions, we rename $u_i = v_i ,\, i = 2, \ldots, n.$ From this,  it is then clear that the first $m-1$ controls $v_i$ are the Dirichlet controls, while the $v_i , i = m + 1, \ldots , n$,  are the Neumann controls.

We are interested in solving the optimal control problem \eqref{OP2}-\eqref{ST2}.
We let
\begin{align}\label{mp2}
\mathbb{J}(v):=\frac{1}{2}\sum_{i=1}^{n}\int_{Q_i} \abs{y^i-y_d^i}^2\,dxdt+\frac{1}{2}\sum_{i=2}^{n}\int_0^T|v_i|^{2} \,dt,
\end{align}
where $y=(y^i)_i$ satisfies \eqref{ST2} and $y_d=(y_d^i)_i\in L^2((0,T);\mathbb{L}^2)$.  Recall that $\mathbb{U}_{ad}$ is a closed and convex subset of $ \left(L^2(0,T)\right)^{n-1}.$
We also let
$$
\mathbb{H}_{a}^\alpha:=\dis \prod_{i= 1}^n H_{a^+}^\alpha (a,b_i).$$
We endow $\mathbb{L}^2$ and $\mathbb{H}_{a}^\alpha$  with the norms given, respectively, by
\begin{equation}\label{normHstar}
\|\rho\|^2_{\mathbb{L}^2}=\dis \sum_{i=1}^n\|\rho^i\|^2_{L^2(a,b_i)},\quad \rho=(\rho^i)_i\in \mathbb{H},
\end{equation}
and
\begin{equation}\label{normHastar}
\|\rho\|^2_{\mathbb{H}_{a}^\alpha}=\dis \sum_{i=1}^n\|\rho^i\|^2_{H^\alpha_{a^+}(a,b_i)}=\sum_{i=1}^n\left(
\norm{\rho^i}^2_{L^2(a,b_i)}+\norm{\D^{\alpha}_{a^+}\rho^i}^2_{L^2(a,b_i)}\right),\quad \rho=(\rho^i)_i\in \mathbb{H}_{a}^\alpha.\end{equation}

As for the case of a single edge studied in the previous section,  to investigate the minimization problem \eqref{OP2}-\eqref{ST2}, we need some preliminary results.  We start with some existence and regularity results.

\subsection{Homogeneous boundary fractional  Sturm–Liouville parabolic equations in a star graph}	

We  consider the following fractional Sturm-Liouville boundary value problem on a general star graph:
\begin{equation}\label{pa7ad}
\left\{
\begin{array}{lllllllllllllllll}
\displaystyle	\rho_t^i+\mathcal{D}_{b_i^-}^\alpha(\beta^i\mathbb{D}_{a^+}^\alpha \rho^i)+q^i\rho^i&=&g^i&\text{in}& Q_i,\,i=1,\dots, n,\\
\displaystyle	(I_{a^+}^{1-\alpha}\rho^i)(a,\cdot)-(I_{a^+}^{1-\alpha}\rho^j)(a,\cdot)&=&0&\text{in}&(0,T),~i\neq j=1,\dots, n,\\
\displaystyle	\sum_{i=1}^n (\beta^i\mathbb{D}_{a^+}^\alpha \rho^i)(a,\cdot)&=&0&\text{in}& (0,T), \\
\displaystyle	(I_{a^+}^{1-\alpha} \rho^i)(b_i^-,\cdot)&=&0&\text{in}& (0,T), \; i=1,\dots, m\\
\displaystyle	(\beta^i\mathbb{D}_{a^+}^\alpha \rho^i)(b_i^-,\cdot)&=&0&\text{in}&(0,T), \; i=m+1,\dots,n,\\
\displaystyle	\rho^i(\cdot,0) &=&\rho^{0,i}&\text{in}& (a,b_i),~~i=1,\dots,n,
	\end{array}
	\right.
\end{equation}
where $g=(g^i)_i\in L^2((0,T), \mathbb{L}^2)$,  and $\rho^0=(\rho^{0,i})_i\in \mathbb{L}^2$.  

We let
\begin{align}\label{DefVstar}
\mathbb{V}:=\dis \Big\{\rho:=(\rho^i)_i\in\mathbb H_a^\alpha: \;  (I^{1-\alpha}_{a^+}\rho^i)(a)-(I^{1-\alpha}_{a^+}\rho^j)(a)=0, i\neq j,\,\,i,j=1,\ldots, n, \notag\\
\qquad \qquad\qquad\mbox{ and }\; (I_{a^+}^{1-\alpha} \rho^i)(b_i^-)= 0, \; i=1,\dots,m\Big\},
\end{align}
and we endow it  with the norm given by
\begin{equation}\label{normstar}
\| \rho\|^2_{\mathbb{V}}:=\sum_{i=1}^n\|\rho^i\|^2_{H^\alpha_{a^+}(a,b_i)}=\|\rho\|_{\mathbb H_a^\alpha}^2.
\end{equation}

\begin{remark}
We notice that the norm in  $L^2((0,T); \mathbb{V})$ is given by
\begin{align}\label{normL2DA}
\|\rho\|_{L^2((0,T); \mathbb{V})}=\dis \int_0^T \|\rho(\cdot,t)\|^2_{\mathbb{V}}\;dt
=\dis \sum_{i=1}^n \int_0^T \left(\|\rho^i(\cdot,t)\|_{L^2(a,b_i)}^2+ \|\Dr_{a^+}^\alpha \rho^i(\cdot,t)\|_{L^2(a,b_i)}^2\right)\;dt.
\end{align}
\end{remark}

We make the following assumption.

\begin{assumption}\label{assump2}
We assume that $q^i\in L^\infty(a,b_i)$,  $\beta^i\in C([a,b_i])$ ($i=1,\ldots,n$),  and there are constants $\beta^{0,i}$ and $q^{0,i}$ such that
\begin{align*}
\beta^i(x)\ge \beta^{0,i}> 0\;\mbox{ for all } x\in [a,b_i]\,\mbox{ and }\; q^i(x)\geq q^{0,i}> 0 \;\mbox{ for a.e. }\; x\in (a,b_i),\,\, i=1,\cdots,n.
\end{align*}
\end{assumption}

From now on, we set
\begin{align} 
\underline{q^0}:=\dis \min_{1\leq i\leq n}q^{0,i},\quad
\underline{\beta^0}:=\dis\min_{1\leq i\leq n}\beta^{0,i},\label{notations1}\\
\overline{q}:=\dis \max_{1\leq i\leq n}\|q^i\|_{\infty},\quad
\overline{\beta}:=\dis\max_{1\leq i\leq n}\|\beta^i\|_{\infty},\label{notations2}
\end{align}
where $\|\beta^i\|_{\infty}:=\dis \max_{x\in[a,b_i]}|\beta^i(x)| \hbox{ and } \dis \|q^i\|_{\infty}:= \sup_{x\in(a,b_i)}|q^i(x)|.$

\begin{remark}\label{rem-SG2}
As in Lemma \ref{lembil} and Remark \ref{rem-SG} we have the following.
\begin{enumerate}[(a)]
\item Consider the bilinear and symmetric form $\mathbb E:\mathbb V\times\mathbb V\to \R$ given by
\begin{align*}
\mathbb E(\rho,\phi):= \sum_{i=1}^n\int_{Q_i} \beta^i(x)\Dr_{a^+}^\alpha \phi^i (x,t) \Dr_{a^+}^\alpha \rho^i(x,t)\;dx dt
+  \sum_{i=1}^n\int_{Q_i}q^i(x)\phi^i(x,t)\rho^i(x,t)\;dxdt.
\end{align*}
It is straightforward to show that the form $\mathbb E$ is continuous and coercive. In addition, it is closed in $\mathbb L^2$.  Let $\mathbb A$ be the selfadjoint operator on $\mathbb L^2$ associated with $\mathbb E$ in the sense that
\begin{equation*}
\begin{cases}
D(\mathbb A)=\{\rho\in\mathbb V:\; \exists\; F\in\mathbb L^2 \mbox{ such that } \mathbb E(\rho,\phi)=(F,\phi)_{\mathbb L^2}\;\forall\;\phi\in\mathbb V\}\\
\mathbb A\rho=F.
\end{cases}
\end{equation*}
Then,  the operator $-\mathbb A$ generates a strongly continuous and analytic semigroup $\mathbb S=(\mathbb S(t))_{t\ge 0}$ in $\mathbb L^2$.
\item The operator $\mathbb A$ can be also viewed as a bounded operator from $\mathbb V$ into $\mathbb V^\star$ given by
\begin{align*}
\langle \mathbb A\rho,\phi\rangle_{\mathbb V^\star,\mathbb V}=\mathbb E(\rho,\phi) \;\mbox{ with }  \mathbb A\rho=(\mathbb A\rho)_i^i \mbox{ where } (\mathbb A\rho)^i=\mathcal{D}_{b_i^-}^\alpha(\beta^i\mathbb{D}_{a^+}^\alpha \rho^i)+q^i\rho^i.
\end{align*}
In this sense,  the operator $-\mathbb A$ also generates a strongly continuous and analytic semigroup in $\mathbb V^\star$. As before,  if there is no confusion we shall use the same notation for the above defined operators and semigroups.
\item The system \eqref{pa7ad} can be also rewritten as the following abstract Cauchy problem:
\begin{equation}\label{ACP-2}
\begin{cases}
 \rho_t+\mathbb A\rho=g \;\;&\mbox{ in }\; \dis\prod_{i=1}^nQ_i,\\
\rho(\cdot,0)=\rho^0 &\mbox{ in }\dis \prod_{i=1}^n(a,b_i).
\end{cases}
\end{equation}
\end{enumerate}
\end{remark}

Now, we introduce our notion of solutions to the system \eqref{pa7ad} (and hence to the Cauchy problem \eqref{ACP-2}).

\begin{definition}\label{weaksolutionhomostar}
Let $0<\alpha\leq 1$,   $g=(g^i)_i\in  L^2((0,T), \mathbb{L}^2)$,  and $\rho^0=(\rho^{0,i})_i\in \mathbb{L}^2$ be given.  A function
$
\rho=(\rho^i)_i\in L^2((0,T); \mathbb{V})\cap H^1((0,T);\mathbb V^\star)$
is said to be a weak solution of \eqref{pa7ad} if the  equality
\begin{equation}\label{varadj1}
\mathbb{F}(\rho,\phi)=\mathbb{L}(\phi),\quad \forall \phi\in \mathbb{H}(Q),
\end{equation}
holds, where
$$
\mathbb{H}(Q)=\dis \Big\{\phi=(\phi^i)_i:\;\phi\in L^2((0,T); \mathbb{V})\cap H^1((0,T);\mathbb L^2): 
\;  \phi^i(\cdot,T)=0 ~~\text{in}~~ (a,b_i),\, i=1,\cdots,n \Big\},
$$
 $\mathbb F:L^2((0,T); \mathbb{V})\times \mathbb H(Q)\to \mathbb R$ is given by
\begin{align*}
\mathbb{F}(\rho,\phi)=&-\dis \sum_{i=1}^n\int_{Q_i}\rho^i (x,t)\phi_t^i(x,t)\;dx dt
+
 \sum_{i=1}^n\int_{Q_i} \beta^i(x)(\Dr_{a^+}^\alpha \phi^i )(x,t) (\Dr_{a^+}^\alpha \rho^i)(x,t)\;dx dt\\
&+  \sum_{i=1}^n\int_{Q_i}q^i(x)\phi^i(x,t)\rho^i(x,t)\;dxdt,
\end{align*}
and $\mathbb L:\mathbb H(Q)\to \R$ is defined by
$$
\mathbb{L}(\phi)= \dis \sum_{i=1}^n\int_{Q_i} g^i(x,t) \phi^i(x,t)\;dxdt
+\dis \sum_{i=1}^n\int_a^{b_i}\rho^{0,i}(x)\phi^i(x,0)\;dx.
$$
\end{definition}

We have the following existence result.

\begin{theorem}\label{theoexisthomostar}
 Let $0<\alpha\leq 1$,   $g=(g^i)_i\in  L^2((0,T), \mathbb{L}^2),$ and $\rho^0=(\rho^{0,i})_i\in \mathbb{L}^2$.  Let $q^i$ and $\beta^i$  satisfy Assumption \ref{assump2}.   Then,  there exists a unique weak solution $\rho=(\rho^i)_i\in L^2((0,T);\mathbb{V})\cap  H^1((0,T), \mathbb{V}^\star)$  of \eqref{pa7ad} in the sense of Definition \ref{weaksolutionhomostar}. Moreover,  there is a constant $C=C(q,\beta)>0$ such that the following estimates hold:
\begin{align}
\dis \|\rho(\cdot,T)\|^2_{\mathbb{L}^2}&\leq C\left(
 \|\rho^0\|^2_{\mathbb{L}^2}+\|g\|^2_{L^2((0,T);\mathbb{L}^2)}\right)
 \label{estimation1adjointstateT}\\
 \|\rho\|^2_{L^2((0,T);\mathbb{V})}&\leq C\left(\|\rho^0\|^2_{\mathbb{L}^2}+
 \|g\|^2_{L^2((0,T);\mathbb{L}^2)}\right)\label{estimation1adjointstaterho}\\
   \|\rho_t\|^2_{L^2((0,T);\mathbb{V}^\star)}&\leq C\left(\|\rho^0\|^2_{\mathbb{L}^2}+
 \|g\|^2_{L^2((0,T);\mathbb{L}^2)}\right).\label{mmww}
\end{align}
\end{theorem}

\begin{proof} We proceed in four steps.

\textbf{Step 1.} First, we prove existence. We endow the pre-Hilbert $\mathbb{H}(Q)$ with the norm given by
$$
\||\varphi\||^2=\|\varphi\|^2_{L^2((0,T),\mathbb{V})}+
\|\varphi(\cdot,0)\|^2_{\mathbb{L}^2} ,\;\;\;\forall\; \varphi\in \mathbb{H}(Q).
$$
Then,  $\|\varphi\|_{L^2((0,T),\mathbb{V})}\leq \||\varphi\||$  for every $\varphi\in  \mathbb{H}(Q)$.
That is,  the embedding $\mathbb H(Q)\hookrightarrow L^2((0,T),\mathbb{V})$ is continuous.

Using Cauchy-Schwarz's  inequality,  and the norm on $\Ha$ given by \eqref{normHa}, we get that there is a constant $C>0$ such that
for every $(\rho,\phi)\in L^2((0,T),\mathbb{V})\times\mathbb H(Q)$
\begin{align}\label{varadj1F0}
\left|\mathbb{F}(\rho,\phi)\right|\leq& \dis \sum_{i=1}^n\int_0^T\|\rho^i(\cdot,t) \|_{L^2(a,b_i)}\|\phi_t^i(\cdot,t)\|_{L^2(a,b_i)}\;dt\notag\\
&+
\dis \sum_{i=1}^n\int_0^T\|\beta^i\|_\infty \|\Dr_{a^+}^\alpha \phi^i (\cdot,t)\|_{L^2(a,b_i)}\| \Dr_{a^+}^\alpha \rho^i(\cdot,t)  \|_{L^2(a,b_i)} \;dt\notag\\
&+ \dis \sum_{i=1}^n\int_0^T\|q^i\|_{\infty}\|\phi^i(\cdot,t)\|_{L^2(a,b_i)}\|\rho^i(\cdot,t) \|_{L^2(a,b_i)}\,dt\notag\\
\leq& \dis \sum_{i=1}^n\int_0^T\|\rho^i(\cdot,t) \|_{H_{a^+}^\alpha(a,b_i)}\|\phi_t^i(,\cdot,t)\|_{L^2(a,b_i)}\;dt\notag\\
&+
\dis \sum_{i=1}^n\int_0^T\|\beta^i\|_\infty \|\Dr_{a^+}^\alpha \phi^i (\cdot,t)\|_{L^2(a,b_i)}\| \Dr_{a^+}^\alpha \rho^i(\cdot,t)  \|_{L^2(a,b_i)}\; dt\notag\\
&+ \dis \sum_{i=1}^n\int_0^T\|q^i\|_{\infty}\|\phi^i(\cdot,t)\|_{L^2(a,b_i)}\|\rho^i(\cdot,t) \|_{L^2(a,b_i)}\,dt\notag\\
\leq& C\left( \|\phi_t\|_{L^2((0,T);\mathbb{L}^2)}+ \|\phi\|_{L^2((0,T);\mathbb{V})}\right) \|\rho \|_{L^2((0,T),\mathbb{V})}.
\end{align}
We have shown that for  every fixed $\phi\in \mathbb{H}(Q)$, the functional  $\rho\mapsto \mathbb{F}(\rho,\phi)$ is continuous  on the Hilbert space $L^2((0,T);\mathbb{V})$.

Next, for every $\phi\in \mathbb{H}(Q),$ we have that
\begin{align*}
\left|\mathbb{F}(\phi,\phi)\right|=&-\dis \sum_{i=1}^n\int_{Q_i}\phi^i (x,t)\phi_t^i(x,t) \;dx dt
+
\dis \sum_{i=1}^n\int_{Q_i} \beta^i(x)\left|\Dr_{a^+}^\alpha \phi^i (x,t) \right|^2 \;dx dt\\
&+ \dis \sum_{i=1}^n\int_0^T\int_a^{b_i}q^i(x)\left|\phi^i(x,t)\right|^2\; dx\,dt\\
\geq &\dis \frac 12 \sum_{i=1}^n\|\phi^i (\cdot,0)\|_{L^2(a,b_i)}^2
+
\dis \underline{q^0}\sum_{i=1}^n\int_{Q_i} \left|\Dr_{a^+}^\alpha \phi^i (x,t) \right|^2 \;dx dt
+ \dis \underline{\beta^0}\sum_{i=1}^n\int_{Q_i}\left|\phi^i(x,t)\right|^2 \;dxdt\\
\geq &\dis \frac 12\|\phi (\cdot,0)\|_{\mathbb{L}^2}^2 + \dis \min\left[\underline{q^0},\underline{\beta^0}\right]
\dis \|\phi\|^2_{L^2((0,T);\mathbb{V})}.
\end{align*}
We have shown  that there is a constant $C>0$ such that
for every $\phi\in \mathbb{H}(Q),$
\begin{equation}\label{varadj1F01}
\left|\mathbb{F}(\phi,\phi)\right|\geq C\||\phi \||,
\end{equation}
and this implies that the bilinear form $\mathbb F$ is coercive on $\mathbb H(Q)$.

In addition, we have that for every $\phi\in \mathbb{H}(Q)$,
\begin{align}\label{varadj1L0}
\left|\mathbb{L}(\phi)\right|\leq & \dis \sum_{i=1}^n\int_0^T\| g^i(\cdot,t)\|_{L^2(a,b_i)} \|\phi^i(\cdot,t)\|_{L^2(a,b_i)} \,dt
+\dis \sum_{i=1}^n\|\rho^{0,i}\|_{L^2(a,b_i)} \|\phi^i(\cdot,0)\|_{L^2(a,b_i)} \notag\\
\leq&\| g\|_{L^2((0,T);\mathbb{L}^2)} \|\phi\|_{L^2((0,T);\mathbb{V})}+
\|\rho^{0}\|_{\mathbb{L}^2} \|\phi(\cdot,0)\|_{\mathbb{L}^2}\notag\\
\leq &\left(\| g\|^2_{L^2((0,T);\mathbb{L}^2)}+\|\rho^{0}\|^2_{\mathbb{L}^2} \right)^{1/2} \||\phi\||.
\end{align}
We have shown that the functional $\mathbb{L}$ is continuous on $\mathbb{H}(Q)$.
In view of \eqref{varadj1F0}, \eqref{varadj1F01}  and \eqref{varadj1L0},  we can deduce from Theorem \ref{Theolions61} again that there exists  a $\rho=(\rho^i)_i\in L^2((0,T);\mathbb{V})$ solution of \eqref{pa7ad} in the sense of Definition \ref{weaksolutionhomostar}.  

\textbf{Step 2.} We show that $ \rho_t=( \rho_t^i)_i \in L^2((0,T);\mathbb V^\star)$. We proceed as in Step 2 in the proof of Theorem \ref{Thm1}.
Notice that by Remark \ref{rem-SG2}, $\mathcal{D}_{b_i^-}^\alpha(\beta^i \mathbb{D}_{a^+}^\alpha \rho^i(\cdot,t))+q^i(x)\rho^i(\cdot,,t)\in (\Hi)^\star$ for a.e. $t\in (0,T)$ and for every $i=1,\ldots,n$. Since by assumption $g^i\in L^2((0,T);L^2(a,b_i))$ for all $i=1,\ldots,n$,  and the system \eqref{pa7ad}  can be rewritten as the abstract Cauchy problem \eqref{ACP-2}, we can also deduce as in Section \ref{oneedge} that $\rho(\cdot,t)\in \mathbb V^\star$ for a.e. $t\in (0,T)$, that is,  $\rho_t^i\in (\Hi)^\star$ for a.e. $t\in (0,T)$ and for every $i=1,\ldots,n$.  In addition, it follows from \cite[Chapter IV, Section 1]{lions1961} that $\rho$ is a weak solution of \eqref{pa7ad} in the sense of Definition \ref{weaksolutionhomostar}, that is, \eqref{varadj1} holds, if and only if the following identity 
 \begin{align}\label{w-sol-muti}
\langle \rho_t(\cdot,t),\Phi\rangle_{\mathbb V^\star,\mathbb V}+\mathbb E(\rho(\cdot,t),\Phi)=\langle g(\cdot,t),\Phi\rangle_{\mathbb V^\star,\mathbb V}=(g(\cdot,t),\Phi)_{\mathbb L^2}
\end{align}
holds for every $\Phi\in  \mathbb V$,  and a.e.  $t\in (0,T)$.
  Let $\phi=(\phi^i)_i\in L^2((0,T);\mathbb V)$.
Taking $\Phi:=\phi(\cdot, t)$ as a test function in \eqref{w-sol-muti},
using the continuity of the bilinear form $\mathbb E$, and integrating over $(0,T)$, we obtain that there is a constant $C>0$ (depending only on the coefficients of the operator) such that
\begin{align*}
\left|\int_0^T\langle \rho_t,\phi\rangle_{\mathbb V^\star,\mathbb V}\;dt \right|\leq& C \int_0^T\|\rho(\cdot,t)\|_{\mathbb V}\|\phi(\cdot,t)\|_{\mathbb V}\;dt+\int_0^T\|g(\cdot,t)\|_{\mathbb L^2}\|\phi(\cdot,t)\|_{\mathbb L^2}\;dt\\
\leq& C\left(\|\rho\|_{L^2((0,T);\mathbb V)}+\|g\|_{L^2((0,T);\mathbb L^2)}\right)\|\phi\|_{L^2((0,T);\mathbb V)}.
\end{align*}
This shows that $\rho_t\in L^2((0,T);\mathbb V^\star)$ and 
\begin{equation}\label{mmww2}
\|\rho_t\|_{ L^2((0,T);\mathbb V^\star)}\leq C( \|\rho\|_{L^2((0,T);\mathbb V)}+\|g\|_{L^2((0,T);\mathbb V)}).
\end{equation} 

\textbf{Step 3.} We show uniqueness.
Assume that there exist $\rho_1=(\rho_1^i)_i$ and $\rho_2=(\rho_2^i)_i$ two weak solutions of \eqref{pa7ad} with the same given data. Set $z^i:=\rho_1^i-\rho_2^i,\, 1\leq i\leq n.$ Then, $z=(z^i)_i$ satisfies
\begin{equation}\label{pa7adzero}
\left\{
\begin{array}{lllllllllllllllll}
\displaystyle	z_t^i+\mathcal{D}_{b_i^-}^\alpha(\beta^i\mathbb{D}_{a^+}^\alpha z^i)+q^iz^i&=&0&\text{in}& Q_i,\,i=1,\dots, n,\\
\displaystyle	(I_{a^+}^{1-\alpha}z^i)(a,\cdot)-(I_{a^+}^{1-\alpha}z^j)(a,\cdot)&=&0&\text{in}&(0,T),~i\neq j=1,\dots, n,\\
\displaystyle	\sum_{i=1}^n(\beta^i\mathbb{D}_{a^+}^\alpha z^i)(a,\cdot)&=&0&\text{in}& (0,T), \\
\displaystyle	(I_{a^+}^{1-\alpha} z^i)(b_i^-,\cdot)&=&0&\text{in}& (0,T), \; i=1,\dots, m\\
\displaystyle	(\beta^i\mathbb{D}_{a^+}^\alpha z^i)(b_i^-,\cdot)&=&0&\text{in}&(0,T), \; i=m+1,\dots,n,\\
\displaystyle	z^i(\cdot,0) &=&0&\text{in}& (a,b_i),~~i=1,\dots,n.
	\end{array}
	\right.
\end{equation}
Proceeding as in Step 3 in the proof of Theorem \ref{Thm1} or Step 2 in the proof of Theorem \ref{theoexistnh} by using the characterization \eqref{w-sol-muti}, we can deduce that
$$0\geq \frac{1}{2}\|z\|^2_{\mathbb{L}^2}+\underline{q^0}\|z\|_{L^2((0,T);\mathbb{L}^2)}.$$
Thus, $z=0$ in $\dis \prod_{i=1}^n Q_i$.  We have that  $\rho_1^i=\rho_2^i$ in $Q_i$, $i=1,\cdots,n$, and this implies uniqueness.

\textbf{Step 4.} We prove the estimates \eqref{estimation1adjointstateT}, \eqref{estimation1adjointstaterho}, and \eqref{mmww}.  We proceed as in the proof of Theorems  \ref{Thm1} or Theorem \ref{theoexistnh}. Taking $\Phi:=\rho(\cdot,t)$ as a test function in \eqref{w-sol-muti} and integrating over $(0,T)$,
we get
$$\begin{array}{llll}
\dis \sum_{i=1}^n\frac{1}{2}\left(\|\rho^i(\cdot,T)\|^2_{L^2(a,b_i)}-\|\rho^{0,i}\|^2_{L^2(a,b_i)}\right)+
\dis \sum_{i=1}^n\int_{Q_i} \beta^i(x)\left|\Dr_{a^+}^\alpha \rho^i )(x,t) \right|^2 \;dx dt\\
+ \dis \sum_{i=1}^n\int_{Q_i}q^i(x)\left|\rho^i(x,t)\right|^2 \;dxdt=
 \dis \sum_{i=1}^n\int_{Q_i} g^i(x,t)\rho^i(x,t)\;dxdt,
 \end{array}
 $$
 from which we can deduce by using Young's inequality that
 \begin{align*}
 \dis \frac{1}{2}\|\rho(\cdot,T)\|^2_{\mathbb{L}^2}+
 \dis \min\left(\underline{q^0},\underline{\beta^0}\right)
\dis \|\rho\|^2_{L^2((0,T);\mathbb{V})}
\leq
\dis \frac{1}{2}\|\rho^0\|^2_{\mathbb L^2}+\dis \frac{1}{2\delta}\|g\|^2_{L^2((0,T);\mathbb{L}^2)}+
\frac{\delta}{2}\|\rho\|^2_{L^2((0,T);\mathbb{V})},
\end{align*}
for every $\delta>0$. Letting $\delta:=\min\left(\underline{q^0},\underline{\beta^0}\right)$ we get the estimates \eqref{estimation1adjointstateT}-\eqref{estimation1adjointstaterho}.  The last estimate \eqref{mmww} is obtained by using \eqref{estimation1adjointstateT}-\eqref{estimation1adjointstaterho} and \eqref{mmww2}.
The proof is finished.
\end{proof}

\begin{remark}
Here also, we observe that the unique weak (strong)  solution $\rho$ of the \eqref{pa7ad} is given by
\begin{align*}
\rho(\cdot,t)=\mathbb S(t)\rho^0+\int_0^t\mathbb S(t-s)g(\cdot,s)\;ds,
\end{align*}
where $\mathbb S$ is the strongly continuous and analytic semigroup on $\mathbb L^2$ generated by the operator $\mathbb A$ introduced in Remark \ref{rem-SG2}.
\end{remark}

We have the following regularity result.

\begin{proposition}\label{regular}
Let $0<\alpha\leq 1$ and $g=(g^i)_i\in  L^2((0,T), \mathbb{L}^2)$. Let $\mathbb A$ be the operator defined in Remark \ref{rem-SG2}.
If $\rho^0=(\rho^{0,i})_i\in D(\mathbb A)$, then the unique weak solution $\rho=(\rho^i)_i$ of the system \eqref{pa7ad} is a strong solution, that is, it belongs to $ C([0,T];D(\mathbb{A}))\cap  H^1((0,T), \mathbb{L}^2)$, and the first equation in  \eqref{pa7ad}  is satisfied pointwise.  In particular, we have that $\mathcal{D}_{b_i^-}^\alpha(\beta^i\mathbb{D}_{a^+}^\alpha \rho^i)+q^i\rho^i \in L^2(Q_i) $ for all $i=1,\ldots,n$.
 \end{proposition}

\begin{proof}
We have shown in Remark \ref{rem-SG2} that $\rho$ is a solution of the abstract Cauchy problem \eqref{ACP-2}.  As we have observed  in Remark \ref{existy0T}(d), if $\rho^0\in D(\mathbb A)$,  then $\rho$ is a strong solution in the sense that  $\rho\in C([0,T];D(\mathbb A))\cap  H^1((0,T), \mathbb{L}^2)$ and the first equation in \eqref{pa7ad} is satisfied pointwise.   In particular we have that  $\rho_t\in   L^2((0,T), \mathbb{L}^2)$.
Since $g, \rho_t\in   L^2((0,T), \mathbb{L}^2)$,  it follows from the first equation in \eqref{pa7ad} that $\mathcal{D}_{b_i^-}^\alpha(\beta^i\mathbb{D}_{a^+}^\alpha \rho^i)+q^i\rho^i \in L^2(Q_i)$ for all $i=1,\ldots,n$. The proof is finished.
\end{proof}

We conclude this section with the following additional regularity result.

\begin{lemma}\label{lemmaajout} 
Let $0<\alpha\leq 1$ and $g=(g^i)_i\in  L^2((0,T), \mathbb{L}^2)$. 
Let $\rho$ be unique strong solution solution of \eqref{pa7ad}  with $\rho^{0}=0$. Then there is a constant $C=C(a,b_i,\alpha,\underline{\beta^0},\underline{q^0},\overline{q},\overline{\beta})>0$ such that
\begin{equation}\label{estdab}
  \sum_{i=1}^n \|\mathbb{D}^\alpha_{a^+}\rho(b_i,\cdot)\|^2_{L^2(0,T)}\leq C\|g\|_{L^2((0,T);\mathbb{L}^2)}^2.
\end{equation}
\end{lemma}

\begin{proof}  
Let $\rho^{0}=0\in D(\mathbb A)$ and $\rho\in  C([0,T];D(\mathbb{A}))\cap  H^1((0,T), \mathbb{L}^2)$ be the unique strong  solution of \eqref{pa7ad}. As we have observed in the proof of Proposition \ref{regular}, we have that 
\begin{equation}\label{octobre1}
  \displaystyle	\rho_t^i+\mathcal{D}_{b_i^-}^\alpha(\beta^i\mathbb{D}_{a^+}^\alpha \rho^i)+q^i\rho^i=g^i\text{ in } Q_i,\,i=1,\dots, n,
\end{equation}
where $g^i,  \rho_t^i, \mathcal{D}_{b_i^-}^\alpha(\beta^i\mathbb{D}_{a^+}^\alpha \rho^i)+q^i\rho^i\in L^2(Q_i), \,\,i=1,\dots, n.$ \\
Observing that $$I^\alpha_{b_i^-}\left(\mathcal{D}_{b_i^-}^\alpha(\beta^i\mathbb{D}_{a^+}^\alpha \rho^i)\right)(x,t)=\beta^i(b_i)\dis \mathbb{D}_{a^+}^\alpha \rho^i(b_i^-,t)-(\beta^i\mathbb{D}_{a^+}^\alpha \rho^i)(x,t),$$
and applying $I^\alpha_{b_i^-}$ to both sides of \eqref{octobre1}, we can deduce  that
$$\dis \mathbb{D}_{a^+}^\alpha \rho^i(b_i^-,t)= I^\alpha_{b_i^-}g^i(x,t)-I^\alpha_{b_i^-}\left(\rho_t^i\right)(x,t)
-I^\alpha_{b_i^-}\left(q^i\rho^i\right)(x,t)+(\beta^i\mathbb{D}_{a^+}^\alpha \rho^i)(x,t).$$
Therefore,
$$\begin{array}{lll}
\dis\sum_{i=1}^n\int_a^{b_i}\int_0^T\left|\mathbb{D}_{a^+}^\alpha \rho^i(b_i^-,t)\right|^2 \;dxdt &\leq & \dis\sum_{i=1}^n\int_a^{b_i}\int_0^T\left|I^\alpha_{b_i^-}g^i(x,t)\right|^2\; dxdt\\
&+&\dis\sum_{i=1}^n\int_a^{b_i}\int_0^T\left|I^\alpha_{b_i^-}\left(\rho_t^i\right)(x,t)
\right|^2 \;dxdt \\
&+&\dis\sum_{i=1}^n\int_a^{b_i}\int_0^T\left|I^\alpha_{b_i^-}\left(q^i\rho^i\right)(x,t)\right|^2 \;dx dt \\
&+&\dis\sum_{i=1}^n\int_a^{b_i}\int_0^T\left|(\beta^i\mathbb{D}_{a^+}^\alpha \rho^i)(x,t)\right|^2 \;dxdt 
\end{array}.$$
This means that 
$$\begin{array}{lll}
\dis \sum_{i=1}^n\int_a^{b_i}\int_0^T\left|\mathbb{D}_{a^+}^\alpha \rho^i(b_i^-,t)\right|^2 dx \, dt 
&\leq & \dis \dis\sum_{i=1}^n\|I^\alpha_{b_i^-}g^i\|_{L^2((0,T); L^2(a,b_i))}^2 \\
&+&\dis \dis\sum_{i=1}^n\|I^\alpha_{b_i^-}\left(\rho_t^i\right)\|_{L^2((0,T); L^2(a,b_i))}^2\\
&+&\dis \dis\sum_{i=1}^n\|I^\alpha_{b_i^-}\left(q^i\rho^i\right)\|_{L^2((0,T); L^2(a,b_i))}^2\\
&+&\dis\sum_{i=1}^n\|(\beta^i\mathbb{D}_{a^+}^\alpha \rho^i)\|_{L^2((0,T); L^2(a,b_i))}^2.
\end{array}
$$
Hence, using Lemma \ref{l1}, we obtain that there is a constant $C>0$ such that
$$\begin{array}{lll}
\dis \dis\sum_{i=1}^n\int_a^{b_i}\int_0^T\left|\mathbb{D}_{a^+}^\alpha \rho^i(b_i^-,t)\right|^2 \;dx dt
&\leq  & \dis C\dis\sum_{i=1}^n\|g^i\|_{L^2((0,T); L^2(a,b_i))}^2 \\
&+&\dis C\sum_{i=1}^n\|\rho_t^i\|_{L^2((0,T); L^2(a,b_i))}^2\\
&+&\dis C\overline{q} \dis\sum_{i=1}^n\|\rho^i\|_{L^2((0,T); L^2(a,b_i))}^2\\
&+&\dis C\dis\sum_{i=1}^n\|(\mathbb{D}_{a^+}^\alpha \rho^i)\|_{L^2((0,T); L^2(a,b_i))}^2,
\end{array}
$$ 
which according to the  definition of the norm in $\mathbb{V}$ gives
\begin{align*}
\dis \left(\dis\sum_{i=1}^n\int_a^{b_i}\int_0^T\left|\mathbb{D}_{a^+}^\alpha \rho^i(b_i^-,t)\right|^2 dx \, dt \right)^{ 1/2}
\leq   \dis C\left(\|g\|_{L^2((0,T); \mathbb{L}^2)}+\|\rho_t\|_{L^2((0,T); \mathbb{L}^2)}+ \|\rho\|_{L^2((0,T); \mathbb{V})}\right).
\end{align*}
Using the previous estimates and \eqref{estimation1adjointstaterho}, we can deduce that 
$$
\dis\sum_{i=1}^n\int_0^T\left|(\mathbb{D}_{a^+}^\alpha \rho^i)(b_i^-,t)\right|^2\; dt 
\leq   \dis C(\|g\|_{L^2((0,T); \mathbb{L}^2)}^2.
$$
We have shown \eqref{estdab} and the proof is finished.
\end{proof}

\subsection{Non-homogeneous fractional Sturm–Liouville parabolic equations  in a star graph}
We now consider the  Sturm-Liouville problem with boundary controls in a general star graph, that is the system \eqref{ST2}. 

Let us give the motivation how we shall introduce our notion of solutions in that case.  For any $f=(f^i)_i\in L^2((0,T);\mathbb{L}^2)$, we consider first the following dual system:
\begin{equation}\label{pa7adf}
\left\{
\begin{array}{lllllllllllllllll}
\displaystyle	-\phi_t^i+\mathcal{D}_{b_i^-}^\alpha(\beta^i\mathbb{D}_{a^+}^\alpha \phi^i)+q^i\phi^i&=&f^i&\text{in}& Q_i,\,i=1,\dots, n,\\
\displaystyle	(I_{a^+}^{1-\alpha}\phi^i)(a,\cdot)-(I_{a^+}^{1-\alpha}\phi^j)(a,\cdot)&=&0&\text{in}&(0,T),~i\neq j=1,\dots, n,\\
\displaystyle	\sum_{i=1}^n (\beta^i\mathbb{D}_{a^+}^\alpha \phi^i)(a,\cdot)&=&0&\text{in}& (0,T), \\
\displaystyle	(I_{a^+}^{1-\alpha} \phi^i)(b_i^-,\cdot)&=&0&\text{in}& (0,T), \; i=1,\dots, m\\
\displaystyle	(\beta^i\mathbb{D}_{a^+}^\alpha \phi^i)(b_i^-,\cdot)&=&0&\text{in}&(0,T), \; i=m+1,\dots,n,\\
\displaystyle	\phi^i(\cdot,T) &=&0&\text{in}& (a,b_i),~~i=1,\dots,n.
	\end{array}
	\right.
\end{equation}

Let $\mathbb V$ be the Hilbert spaces defined in  \eqref{DefVstar} and 
\begin{align}\label{sphi}
\Phi:=\Big\{\varphi:=(\varphi^i)_i\in L^2((0,T);\mathbb{ V})\cap H^1((0,T);\mathbb{L}^2):\;\mathcal{D}_{b_i^-}^\alpha(\beta^i\mathbb{D}_{a^+}^\alpha \varphi^i) +q^i\varphi^i\in L^2(Q_i), \;i=1,\ldots,n\notag\\
\varphi^i(\cdot,T)=0 \mbox{ in } (a,b_i),\; i=1,\ldots,n\Big\}.
\end{align}
We endow $\Phi$ with the norm given by
\begin{equation}\label{normPhi}
  \||\varphi\||^2:=\dis \sum_{i=1}^n \|\mathbb{D}^\alpha_{a^+}\varphi(b_i,\cdot)\|^2_{L^2(0,T)}+ \|\varphi(\cdot,0)\|^2_{\mathbb{L}^2}+ \|\varphi\|_{L^2((0,T);\mathbb V)}^2,\;\;\varphi\in\Phi.
\end{equation}

\begin{remark}
We observe that by changing $t\to T-t$,  we have that solutions of \eqref{pa7adf} enjoin the same regularities as the ones stated in 
 Proposition \ref{regular} and Lemma \ref{lemmaajout}, that is, weak and strong solutions coincide. Thus,  if $\phi=(\phi^i)_{i=1^n}$ is the weak (strong) solution of \eqref{pa7adf}, then $\phi\in \Phi$.
\end{remark}

If we multiply the first equation in \eqref{ST2} by $\phi=(\phi^i)_{i=1}^n$ a strong solution of \eqref{pa7adf},  integrate  by parts over $Q$ (by using the integration by parts formulas given in Section \ref{prelim},  and using the initial and boundary conditions) we get the following identity:
\begin{align}\label{defveryweaksol}
& \sum_{i=1}^n\int_{Q_i}y^i \left(-\phi_t^i +\mathcal{D}_{b_i^-}^\alpha(\beta^i\mathbb{D}_{a^+}^\alpha \phi^i)+q^i\phi^i\right)\;dxdt \notag\\
=&\dis \sum_{i=1}^n\int_{Q_i} f^i \phi^i dx\,dt+\dis \sum_{i=1}^n\int_a^{b_i}y^{0,i}(x)\phi^i(x,0)\; dx\notag\\
&+\dis \sum_{i=m+1}^n\int_0^Tv_i(t)I_{a^+}^{1-\alpha}(\phi^i)(b_i^-,t)\; dt-\dis \sum_{i=2}^m\int_0^Tu_i(t)\mathbb{D}_{a^+}^\alpha (\phi^i)(b_i^-,t) \;dt.
\end{align}
Conversely, if $y$ satisfies \eqref{defveryweaksol} for every $\phi\in\Phi$, then it is also easy to see that \eqref{ST2} holds.  Besides, we notice that the identity \eqref{defveryweaksol} makes sense for every $y=(y^i)_i\in L^2((0,T);\mathbb L^2)$ and $y^0=(y^{0,i})_i\in\mathbb L^2$.

With this motivation, here is our new notion of solutions to the system \eqref{ST2}.
	
\begin{definition}\label{weaksolutionnonhomostar}
Let $0<\alpha\leq 1,$  $f=(f^i)_i\in  L^2((0,T), \mathbb{L}^2),$  and $y^0=(y^{0,i})_i\in \mathbb{L}^2$. Let  $q^i\in L^\infty(a,b_i)$ and $\beta^i\in C([a,b_i])$ satisfy Assumption \ref{assump2}.  Let also  $u_i\in L^2(0,T),\, i=2,\ldots, m$   and  $v_i\in L^2(0,T),\, i=m+1,\ldots, n.$ A function $y=(y^i)_i\in L^2((0,T); \mathbb{L}^2)$
is said to be a very-weak solution (or a solution by transposition) of  \eqref{ST2} if the  identity \eqref{defveryweaksol}
holds for every $\phi\in  \Phi$.
\end{definition}

We have the following existence result.
 
\begin{theorem}\label{theoexisthomostarnh1}
 Let $0<\alpha\leq 1,$  $f=(f^i)_i\in  L^2((0,T), \mathbb{L}^2),$ and $y^0=(y^{0,i})_i\in \mathbb{L}^2$.  Let  $q^i\in L^\infty(a,b_i)$ and $\beta^i\in C([a,b_i])$ satisfy Assumption \ref{assump2}.  Let also  $u_i\in L^2(0,T),\, i=2,\ldots, m$   and  $v_i\in L^2(0,T),\, i=m+1,\ldots, n.$   Then,  there exists a unique very-weak solution $y=(y^i)_i\in L^2((0,T);\mathbb{L}^2)$  of \eqref{ST2} in the sense of Definition \ref{weaksolutionnonhomostar}. Moreover,   there is a constant $C=C(a,b_i,\alpha,\underline{\beta^0},\underline{q^0},\overline{q},\overline{\beta})>0$ such that
  \begin{equation}\label{estimation1starh}
\|y\|_{L^2((0,T);\mathbb{L}^2)}^2\leq C \left(\|f\|^2_{L^2((0,T);\mathbb{L}^2)}+\|y^0\|^2_{\mathbb{L}^2}+\|v\|^2_{L^2(0,T)}+
\|u\|^2_{L^2(0,T)}\right).
\end{equation}
\end{theorem}

\begin{proof} 
The proof is inspired from the results contained in \cite[Pages 71-74]{lions1968}. 

Firstly, we consider the mapping
\begin{align*}
\mathcal W:\Phi\to L^2((0,T);\mathbb{L}^2),\; \phi \mapsto \mathcal W\phi: \mbox{ with }
 \int_Qy\mathcal W\phi\;dxdt=\sum_{i=1}^n\int_{Q_i}y^i \left(-\phi_t^i +\mathcal{D}_{b_i^-}^\alpha(\beta^i\mathbb{D}_{a^+}^\alpha \phi^i)+q^i\phi^i\right)\;dxdt,
\end{align*}
for every $y\in L^2((0,T);\mathbb{L}^2)$. It follows from  Proposition \ref{regular} and Lemma \ref{lemmaajout} that the mapping $\mathcal W$ is an isomorphism.

Secondly,  we consider  the linear functional $\mathcal{M}:\Phi\to\R$ defined by
\begin{align*}
\mathcal{M}(\phi)=&\dis \sum_{i=1}^n\int_{Q_i} f^i(x,t) \phi^i(x,t) \;dxdt+\dis \sum_{i=1}^n\int_a^{b_i}y^{0,i}(x)\phi^i(x,0)\; dx\\
&+\dis \sum_{i=m+1}^n\int_0^Tv_i(t)I_{a^+}^{1-\alpha}(\phi^i)(b_i^-,t)\; dt-\dis \sum_{i=2}^m\int_0^Tu_i(t)\mathbb{D}_{a^+}^\alpha (\phi^i)(b_i^-,t) \;dt.
\end{align*}
Calculating and using the estimates in Lemma \ref{trace}, we get that there is a constant $C>0$ such that
\begin{align}\label{exo6}
\left|\mathcal{M}(\phi)\right|\leq & \|f\|_{L^2((0,T);\mathbb{L}^2)}\|\phi\|_{L^2((0,T);\mathbb{L}^2)}+\dis \|y^0\|_{\mathbb{L}^2}\|\phi(\cdot,0)\|_{\mathbb{L}^2}\notag\\
&+\dis \|v\|_{L^2(0,T)}\sum_{i=1}^n\|I_{a^+}^{1-\alpha}(\phi^i)(b_i^-,\cdot)\|_{L^2(0,T)}
+\dis \dis \|u\|_{L^2(0,T)}\left(\sum_{i=1}^n\|\mathbb{D}_{a^+}^\alpha (\phi^i)(b_i^-,.)\|^2_{L^2(0,T)}\right)^{1/2}\notag\\
\leq & \|f\|_{L^2((0,T);\mathbb{L}^2)}\|\phi\|_{L^2((0,T);\mathbb{V})}+\dis \|y^0\|_{\mathbb{L}^2}\|\phi(\cdot,0)\|_{\mathbb{L}^2}\notag\\
&+C \|v\|_{L^2(0,T)}\|\phi\|_{L^2((0,T);\mathbb{V})}
+\dis \dis \|u\|_{L^2(0,T)}\left(\sum_{i=1}^n\|\mathbb{D}_{a^+}^\alpha (\phi^i)(b_i^-,\cdot)\|^2_{L^2(0,T)}\right)^{1/2}\notag\\
\leq&C \left(\|f\|^2_{L^2((0,T);\mathbb{L}^2)}+\|y^0\|^2_{\mathbb{L}^2}+\|v\|^2_{L^2(0,T)}+\|u\|^2_{L^2(0,T)}
\right)^{1/2}\times\notag \\
&\left(\|\phi\|^2_{L^2((0,T);\mathbb{V})}+ \|\phi(\cdot,0)\|^2_{\mathbb{L}^2}+\sum_{i=1}^n\|\mathbb{D}_{a^+}^\alpha (\phi^i)(b_i^-,\cdot)\|^2_{L^2(0,T)}\right)^{1/2}\notag\\
\leq& C\left(\|f\|^2_{L^2((0,T);\mathbb{L}^2)}+\|y^0\|^2_{\mathbb{L}^2}+\|v\|^2_{L^2(0,T)}+\|u\|^2_{L^2(0,T)}
\right)^{1/2}\Big)\||\phi\||.
\end{align}
We have shown that the linear functional $\mathcal{M}$ is  continuous on $\Phi$.  Consequently,
it follows from the results contained in \cite[Pages 71-74]{lions1968} that
there exists a unique function $y=(y^i)_i\in L^2((0,T); \mathbb{L}^2)$ such that \eqref{defveryweaksol} holds true for every $\phi=(\phi^i)_i\in \Phi.$ We have shown that the system \eqref{ST2} as a unique very-weak solution $y=(y^i)_i\in L^2((0,T); \mathbb{L}^2)$  in the sense of Definition \ref{weaksolutionnonhomostar}.

Finally,  taking in \eqref{defveryweaksol} $\phi=(\phi^i)_i$  solution of \eqref{pa7adf} with $f=y$, and using \eqref{exo6}, we obtain that there is a constant $C>0$ such that 
\begin{equation}\label{JJ1}
\|y\|^2_{L^2((0,T);\mathbb{L}^2)} \leq C\left(\|f\|^2_{L^2((0,T);\mathbb{L}^2)}+\|y^0\|^2_{\mathbb{L}^2}+\|v\|^2_{L^2(0,T)}+
\|u\|^2_{L^2(0,T)}
\right)^{1/2}\||\phi\||.
\end{equation}
Noticing that there is a constant $C>0$ such that $\||\phi\||\leq C\|f\|_{L^2((0,T);\mathbb L^2)}=C\|y\|_{L^2((0,T);\mathbb L^2)}$ and using
\eqref{JJ1},  we can deduce that there is a constant $C>0$ such that
$$\|y\|_{L^2((0,T);\mathbb{L}^2)}^2 \leq C \left(\|f\|^2_{L^2((0,T);\mathbb{L}^2)}+\|y^0\|^2_{\mathbb{L}^2}+\|v\|^2_{L^2(0,T)}+
\|u\|^2_{L^2(0,T)}\right). $$
 We have shown  \eqref{estimation1starh} and the proof is concluded.
\end{proof}

\subsection{Existence of minimizers and optimality conditions in a general star graph}

We now consider the optimal control problem \eqref{OP2}-\eqref{ST2}. Let $\mathbb J$ be the functional defined in \eqref{mp2}, and recall that
$\mathbb{U}_{ad}$ is a closed  and convex subset of $ \left(L^2(0,T)\right)^{n-1}.$

We have the following existence result of optimal controls.

\begin{theorem}\label{Thmm}
Let $0<\alpha\leq 1$.  Let $q^i$ and $\beta^i$  satisfy Assumption \ref{assump2}.  Then,  there exists a unique solution $\hat u\in \mathbb{U}_{ad}$ of the optimal control problem \eqref{OP2}-\eqref{ST2}.
\end{theorem}

\begin{proof} Notice that $\mathbb{J}(v)\geq  0$ for all $v\in \mathbb{U}_{ad} $ . Let $\{v_k \}=\{(v_{i,k}), i=2,\cdots,n\} \subset \mathbb{U}_{ad}$ be a minimizing sequence such that	
	\begin{align*}
	\underset{k \to \infty}{\lim} \mathbb{J}(v_k ) =\underset{v \in \mathbb{U}_{ad}}{\min} \mathbb{J}(v).
	\end{align*}
	Due to the non-negativity of the functional $\mathbb J$, such a minimizing sequence always exists.
	Then,  there is a constant $C > 0$ (independent of $k$) such that
	\begin{align}\sum_{i=2}^n\label{eq1i}\norm{v_{i,k}}_{L^2(0,T)} \leq C.
	\end{align}
	
	The control $v_{i,k}$ is associated with the state $y_k^i, i=2,\cdots,n,$ with is a very-weak solution of
	\begin{equation}\label{pa2i}
\left\{	
	\begin{array}{lllllll}
	\displaystyle (y_k)_t^i+\mathcal{D}_{b_i^-}^\alpha(\beta^i\mathbb{D}_{a^+}^\alpha y_k^i)+q^iy_k^i&=&f^i\,\,&\text{in}&~ Q_i,\,i=1,\dots, n,\\
\displaystyle	(I_{a^+}^{1-\alpha}y_k^i)(a,\cdot)-(I_{a^+}^{1-\alpha}y_k^j)(a,\cdot)&=&0&\text{in}&~(0,T),~i\neq j=1,\dots, n, \\
\displaystyle	\sum_{i=1}^n(\beta^i\mathbb{D}_{a^+}^\alpha y_k^i)(a,\cdot)&=&0~~~~&\text{in}&~ (0,T), \\
\displaystyle	(I_{a^+}^{1-\alpha} y_k^1)(b_1^-,\cdot)&=&0~~~~&\text{in}&~ (0,T), \\
\displaystyle	(I_{a^+}^{1-\alpha} y_k^i)(b_i^-,\cdot)&=&v_{i,k}~~~&\text{in}&~ (0,T),\; i=2,\dots, m\\
\displaystyle(	\beta^i\mathbb{D}_{a^+}^\alpha y_k^i)(b_i^-,\cdot)&=&v_{i,k}~~~&\text{in}&~(0,T), \; i=m+1,\dots,n,\\
\displaystyle	y_k^i(\cdot,0) &=&y^{0,i}~&\text{in}&~ (a,b_i),~~i=1,\dots,n.
	\end{array}
	\right.
\end{equation}
It follows from \eqref{estimation1starh} and \eqref{eq1i} that there is a constant $C>0$ (independent of $k$) such that
	\begin{align}\label{eq2i}
	\norm{y_k}_{L^2((0,T);\mathbb{L}^2)}\leq C.
	\end{align}
	From \eqref{eq1i} and \eqref{eq2i},  we have that there exists $\hat{u} \in \left(L^2(0,T)\right)^{n-1}$ and $\hat{y} \in L^2((0,T);\mathbb{L}^2)$ such that, as $k\to\infty$, we have
		\begin{align}
		v_k&\rightharpoonup \hu~ \text{weakly in}~ \left(L^2(0,T)\right)^{n-1}\label{eq9i}
		\\
y_k&\rightharpoonup \hy~ \text{weakly in}~L^2((0,T);\mathbb{L}^2)\label{eq10i}.
		\end{align}
	Since $v_k\in \mathbb{U}_{ad}$, which is a closed and convex subset of $\left(L^2(0,T)\right)^{n-1}$, we can deduce that
\begin{equation}\label{ajout10}
  \hu\in \mathbb{U}_{ad}.
\end{equation}
It follows from the definition of very-weak solutions of  \eqref{pa2i} that
\begin{align}\label{defveryweaksoln}
& \sum_{i=1}^n\int_{Q_i}y_k^i \left(-\phi_t^i +\mathcal{D}_{b_i^-}^\alpha(\beta^i\mathbb{D}_{a^+}^\alpha \phi^i)+q^i\phi^i\right)\;dxdt \notag\\
=&\dis \sum_{i=1}^n\int_{Q_i} f^i \phi^i\; dxdt+\dis \sum_{i=1}^n\int_a^{b_i}y^{0,i}(x)\phi^i(x,0) \;dx\notag\\
&-\dis \sum_{i=2}^m\int_0^Tv_{i,k}(t)\mathbb{D}_{a^+}^\alpha (\phi^i)(b_i^-,t)\; dt+\dis \sum_{i=m+1}^n\int_0^Tv_{i,k}(t)(I_{a^+}^{1-\alpha}\phi^i)(b_i^-,t) \;dt,
\end{align}
for every $\phi=(\phi^i)_i\in \Phi$.
Passing to the limit in \eqref{defveryweaksoln} as $k\to\infty$,  while using \eqref{eq9i}-\eqref{eq10i}, we obtain that
$\hy=(\hy^i)_i$ is a very-weak solution of \eqref{pa2i} with $v_{i,k}$ replaced by $\hat u_i$.
The uniqueness follows from the strict convexity of  $\mathbb{J}$.
 The proof is finished.
\end{proof}

Next, we characterize the first order optimality conditions.

\begin{theorem}\label{theo-48}
Let $0<\alpha\leq 1$.  Let $q^i$ and $\beta^i$  satisfy Assumption \ref{assump2}.
	Let $\hu=(\hu_i)_i\in \mathbb{U}_{ad}$ be the optimal control for  the minimization problem  \eqref{OP2}-\eqref{ST2}. Then,  the first order necessary optimality conditions are given by 
\begin{align}\label{pa6i}
& \sum_{i=2}^m\int_0^T \left(\hu_i(t) - (\beta^i\Da \hp^i)(b_i^-,t)\right) ({v}_i(t)-\hu_i(t))\;dt\notag\\
&+
\dis \sum_{i=m+1}^n \int_0^T \left(\hu_i(t) + (I^{1-\alpha}_{a^+}\hp^i)(b_i^-,t)\right)({v}_i(t)-\hu_i(t))\;dt \geq 0
\end{align}
 for all $ {v}=(v_i)_i\in \mathbb{U}_{ad}$, where $\hp$ solves the backward system
 \begin{equation}\label{pa5i}
	\left\{
 \begin{array}{lllllll}
	\displaystyle -\hp_t^i+\mathcal{D}_{b_i^-}^\alpha(\beta^i\mathbb{D}_{a^+}^\alpha \hp^i)+q^i\hp^i&=&\hy^i-y^i_d\,\,&\text{in}&~ Q_i,\,i=1,\dots, n,\\
\displaystyle	(I_{a^+}^{1-\alpha}\hp^i)(a,\cdot)-(I_{a^+}^{1-\alpha}\hp^j)(a,\cdot)&=&0~&\text{in}&~(0,T)~i\neq j=1,\dots, n, \\
\displaystyle	\sum_{i=1}^n(\beta^i\mathbb{D}_{a^+}^\alpha \hp^i)(a,\cdot)&=&0~~~~&\text{in}&~ (0,T), \\
\displaystyle	(I_{a^+}^{1-\alpha} \hp^i)(b_i^-,\cdot)&=&0~~~&\text{in}&~ (0,T), \; i=1,\dots, m\\
\displaystyle	(\beta^i\mathbb{D}_{a^+}^\alpha \hp^i)(b_i^-,\cdot)&=&0~~~&\text{in}&~(0,T), \; i=m+1,\dots,n,\\
\displaystyle	\hp^i(\cdot,T) &=&0~&\text{in}&~ (a,b_i),~~i=1,\dots,n,
	\end{array}
	\right.
\end{equation}
and $\hy$ is the unique very-weak solution of the state equation \eqref{ST2} with $u_i$ and $v_i$ replaced with $\hat u_i$.  In addition, \eqref{pa6i} is equivalent to
\begin{equation}\label{ww2}
\hu=(\hu_i)_i=
\begin{cases}
\mathbb P_{\mathbb U_{ad}}\Big((\beta^i\Dr_{a^+}^\alpha\hp^i)(b_i^-,\cdot)\Big),\quad &i=2,\ldots,m\\
\mathbb P_{\mathbb U_{ad}}\Big(-(I_{a^+}^{1-\alpha}\hp^i)(b_i^-,\cdot)\Big), &i=m+1,\ldots, n,
\end{cases}
\end{equation}
where $\mathbb P_{\mathbb U_{ad}}$ is the projection onto the set $\mathbb U_{ad}$. 
Moreover, since the functional $\mathbb J$ is convex, we have that \eqref{pa6i} is also a sufficient condition.
\end{theorem}

\begin{proof}
The proof follows the lines as the case of a single edge in Theorem \ref{theo-39} by using 
the Lagrangian	
	\begin{align*}
	\mathcal{L}(y,v,p)
=&\mathcal{J}(v)	+\dis \sum_{i=1}^n\int_{Q_i} f^i\,p^i\, dxdt +\sum_{i=1}^n \int_a^{b_i}y^{i0} (x)p^i(x,0)\;dx\\
&-\sum_{i=2}^m \int_0^Tv_i (t) (b_i^-)(\beta^i\D^{\alpha}_{a^+}p^i)(b_i^-,t)\, dt+\dis \sum_{i=m+1}^n \int_0^T  v_i(t)(I^{1-\alpha}_{a^+}p^i)(b_i^-,t)\;dt \\
&-\dis \sum_{i=1}^n\int_{Q_i}\left(-p_t^i +\dis  \Dc_{b^-}^\alpha\,(\beta\, \Dr_{a^+}^\alpha p^i) +q^i p^i\right)y^i \,dxdt.
	\end{align*}
We omit the details for brevity.
\end{proof}

\section{Concluding remarks}
We investigated an optimal control problem of a fractional parabolic partial differential equation involving a fractional Sturm-Liouville operator in a space interval,  and in a general star graph, where the Sturm-Liouville operator is obtained as a composition of a left fractional Caputo derivative,  and a right fractional Riemann–Liouville derivative.  We proved that the considered fractional optimal control in an interval as well as in the graph has a unique solution. We then derived the optimality system that characterizes the control in an edge by means of the Euler-Lagrange optimality conditions,  and also in the graph by using the method of Lagrange multipliers.\\

\noindent
{\bf Acknowledgement:} We would like to thank both referees for their careful reading of the manuscript and their precise comments that helped to improve the final version of the paper.

\bibliographystyle{abbrv}
\bibliography{MW}

\begin{thebibliography}{10}

\bibitem{agrawal2004general}
O.~P. Agrawal.
\newblock A general formulation and solution scheme for fractional optimal
  control problems.
\newblock {\em Nonlinear Dynam.}, 38(1-4):323--337, 2004.

\bibitem{Agr2007}
O.~P. Agrawal.
\newblock Fractional variational calculus in terms of {R}iesz fractional
  derivatives.
\newblock {\em J. Phys. A}, 40(24):6287--6303, 2007.

\bibitem{qasem}
Q.~M. Al-Mdallal.
\newblock On the numerical solution of fractional {S}turm–{L}iouville
  problems.
\newblock {\em International Journal of Computer Mathematics},
  87(12):2837--2845, 2010.

\bibitem{Mehmeti}
F.~Ali~Mehmeti.
\newblock {\em Nonlinear waves in networks}, volume~80 of {\em Mathematical
  Research}.
\newblock Akademie-Verlag, Berlin, 1994.

\bibitem{AGKW-NA}
E.~Alvarez, C.~G. Gal, V.~Keyantuo, and M.~Warma.
\newblock Well-posedness results for a class of semi-linear super-diffusive
  equations.
\newblock {\em Nonlinear Anal.}, 181:24--61, 2019.

\bibitem{Att-BM}
H.~Attouch, G.~Buttazzo, and c.~Michaille.
\newblock {\em Variational analysis in {S}obolev and {BV} spaces}, volume~17 of
  {\em MOS-SIAM Series on Optimization}.
\newblock Society for Industrial and Applied Mathematics (SIAM), Philadelphia,
  PA; Mathematical Optimization Society, Philadelphia, PA, second edition,
  2014.
\newblock Applications to PDEs and optimization.

\bibitem{Baz}
E.~Bajlekova.
\newblock {\em Fractional evolution equations in Banach spaces}.
\newblock Technische Universiteit Eindhoven Eindhoven, 2001.

\bibitem{BK-2013}
G.~Berkolaiko and P.~Kuchment.
\newblock {\em Introduction to quantum graphs}, volume 186 of {\em Mathematical
  Surveys and Monographs}.
\newblock American Mathematical Society, Providence, RI, 2013.

\bibitem{gunter1}
U.~Brauer and G.~Leugering.
\newblock On boundary observability estimates for semidiscretizations of a
  dynamic network of elastic strings.
\newblock {\em Recent advances in control of PDEs. Control \&Cybernetics},
  28(-):421--447, 1999.

\bibitem{schimidt}
Z.~Chen, M.~M. Meerschaert, and E.~Nane.
\newblock On the modelling and exact controllability of networks of vibrating
  strings.
\newblock {\em SIAM Journal of Control and Optimization}, 30(-):229--245, 1992.

\bibitem{Chen}
Z.-Q. Chen, M.~M. Meerschaert, and E.~Nane.
\newblock Space-time fractional diffusion on bounded domains.
\newblock {\em J. Math. Anal. Appl.}, 393(2):479--488, 2012.

\bibitem{dager}
R.~D\'{a}ger.
\newblock Observation and control of vibrations in tree-shaped networks of
  strings.
\newblock {\em SIAM Journal of Control and Optimization}, 43(-):590--623, 2004.

\bibitem{DZ-2006}
R.~D\'{a}ger and E.~Zuazua.
\newblock {\em Wave propagation, observation and control in {$1\text{-}d$}
  flexible multi-structures}, volume~50 of {\em Math\'{e}matiques \&
  Applications (Berlin) [Mathematics \& Applications]}.
\newblock Springer-Verlag, Berlin, 2006.

\bibitem{de-Simon}
L.~de~Simon.
\newblock Un'applicazione della teoria degli integrali singolari allo studio
  delle equazioni differenziali lineari astratte del primo ordine.
\newblock {\em Rend. Sem. Mat. Univ. Padova}, 34:205--223, 1964.

\bibitem{dorville2011optimal}
R.~Dorville, G.~M. Mophou, and V.~S. Valmorin.
\newblock Optimal control of a nonhomogeneous {D}irichlet boundary fractional
  diffusion equation.
\newblock {\em Comput. Math. Appl.}, 62(3):1472--1481, 2011.

\bibitem{Ga-Wa-2021}
C.~G. Gal and M.~Warma.
\newblock {\em Fractional-in-time semilinear parabolic equations and
  applications}, volume~84 of {\em Math\'{e}matiques \& Applications (Berlin)
  [Mathematics \& Applications]}.
\newblock Springer, Cham, [2020] \copyright 2020.

\bibitem{Pod2}
N.~Heymans and I.~Podlubny.
\newblock Physical interpretation of initial conditions for fractional
  differential equations with {R}iemann-{L}iouville fractional derivatives.
\newblock {\em Rheologica Acta}, 45(5):765--771, 2006.

\bibitem{darius}
D.~Idczak and S.~Walczak.
\newblock Fractional {S}obolev spaces via {R}iemann-{L}iouville derivatives.
\newblock {\em J. Funct. Spaces Appl.}, pages Art. ID 128043, 15, 2013.

\bibitem{hassan2015}
H.~Khosravian-Arab, M.~Dehghan, and M.~Eslahchi.
\newblock Fractional {S}turm–{L}iouville boundary value problems in unbounded
  domains: Theory and applications.
\newblock {\em Journal of Computational Physics}, 299:526--560, 2015.

\bibitem{kil}
A.~A. Kilbas, H.~M. Srivastava, and J.~J. Trujillo.
\newblock {\em Theory and applications of fractional differential equations},
  volume 204 of {\em North-Holland Mathematics Studies}.
\newblock Elsevier Science B.V., Amsterdam, 2006.

\bibitem{klimek2013}
M.~Klimek and O.~Agrawal.
\newblock Fractional {S}turm–{L}iouville problem.
\newblock {\em Journal of Mathematical Analysis and Applications},
  66(5):795--812, 2013.

\bibitem{klimek2016}
M.~Klimek, A.~B. Malinowska, and T.~Odzijewicz.
\newblock Applications of the fractional {S}turm-{L}iouville problem to the
  space-time fractional diffusion in a finite domain.
\newblock {\em Fract. Calc. Appl. Anal.}, 19(2):516--550, 2016.

\bibitem{klimek2014}
M.~Klimek, T.~Odzijewicz, and A.~B. Malinowska.
\newblock Variational methods for the fractional {S}turm–{L}iouville problem.
\newblock {\em Journal of Mathematical Analysis and Applications},
  416(1):402--426, 2014.

\bibitem{LG-2004}
J.~E. Lagnese and G.~Leugering.
\newblock {\em Domain decomposition methods in optimal control of partial
  differential equations}, volume 148 of {\em International Series of Numerical
  Mathematics}.
\newblock Birkh\"{a}user Verlag, Basel, 2004.

\bibitem{gunter2}
G.~Leugering.
\newblock Reverberation analysis and control of networks of elastic strings.
\newblock {\em Control of PDE and appl., Lect. Notes in Pure and Applied
  Math.}, -(-):193--206, 1996.

\bibitem{gunter3}
G.~Leugering and G.~Mophou.
\newblock Instantaneous optimal control of friction dominated flow in a
  gas-network.
\newblock In {\em Proceedings of DFG-AIMS Workshop on Shape Optimization,
  Homogenization and Control. AIMS Senegal, Mbour, Senegal}.
  Birk\"{a}user-Verlag, 2017.

\bibitem{lions1961}
J.-L. Lions.
\newblock {\em \'{E}quations diff\'{e}rentielles op\'{e}rationnelles et
  probl\`emes aux limites}.
\newblock Die Grundlehren der mathematischen Wissenschaften, Bd. 111.
  Springer-Verlag, Berlin-G\"{o}ttingen-Heidelberg, 1961.

\bibitem{lions1968}
J.-L. Lions.
\newblock {\em Contr\^{o}le optimal de syst\`emes gouvern\'{e}s par des
  \'{e}quations aux d\'{e}riv\'{e}es partielles}.
\newblock Avant propos de P. Lelong. Dunod, Paris; Gauthier-Villars, Paris,
  1968.

\bibitem{mag}
J.-L. Lions and E.~Magenes.
\newblock {\em Non-homogeneous boundary value problems and applications. {V}ol.
  {II}}.
\newblock Springer-Verlag, New York-Heidelberg, 1972.
\newblock Translated from the French by P. Kenneth, Die Grundlehren der
  mathematischen Wissenschaften, Band 182.

\bibitem{Lumer}
G.~Lumer.
\newblock Connecting of local operators and evolution equations on networks.
\newblock In {\em Potential theory, {C}openhagen 1979 ({P}roc. {C}olloq.,
  {C}openhagen, 1979)}, volume 787 of {\em Lecture Notes in Math.}, pages
  219--234. Springer, Berlin, 1980.

\bibitem{mehra2019}
V.~Mehandiratta, M.~Mehra, and G.~Leugering.
\newblock Existence and uniqueness results for a nonlinear {C}aputo fractional
  boundary value problem on a star graph.
\newblock {\em Journal of Mathematical Analysis and Applications},
  477(2):1243--1264, 2019.

\bibitem{mehra2021}
V.~Mehandiratta, M.~Mehra, and G.~Leugering.
\newblock Existence results and stability analysis for a nonlinear fractional
  boundary value problem on a circular ring with an attached edge: A study of
  fractional calculus on metric graph.
\newblock {\em Networks \& Heterogeneous Media}, 16(2):155--185, 2021.

\bibitem{mehra2021a}
V.~Mehandiratta, M.~Mehra, and G.~Leugering.
\newblock Fractional optimal control problems on a star graph: Optimality
  system and numerical solution.
\newblock {\em Mathematical Control \& Related Fields}, 11(1):189--209, 2021.

\bibitem{MML}
V.~Mehandiratta, M.~Mehra, and G.~Leugering.
\newblock Optimal control problems driven by time-fractional difusion equation
  on metric graphs: optimatily system and finite difference approximation.
\newblock {\em Preprint}, 2021.

\bibitem{mophou2020}
G.~Mophou, G.~Leugering, and P.~S. Fotsing.
\newblock Optimal control of a fractional {S}turm--{L}iouville problem on a
  star graph.
\newblock {\em Optimization}, 70(3):659--687, 2021.

\bibitem{Mophou2011}
G.~M. Mophou.
\newblock Optimal control of fractional diffusion equation.
\newblock {\em Comput. Math. Appl.}, 61(1):68--78, 2011.

\bibitem{JNakagawa}
J.~Nakagawa, K.~Sakamoto, and M.~Yamamoto.
\newblock Overview to mathematical analysis for fractional diffusion
  equations—new mathematical aspects motivated by industrial collaboration.
\newblock {\em Journal of Math-for-Industry}, 2(10):99--108, 2010.

\bibitem{pod}
I.~Podlubny.
\newblock {\em Fractional differential equations}, volume 198 of {\em
  Mathematics in Science and Engineering}.
\newblock Academic Press, Inc., San Diego, CA, 1999.
\newblock An introduction to fractional derivatives, fractional differential
  equations, to methods of their solution and some of their applications.

\bibitem{rivero2013}
M.~Rivero, J.~Trujillo, and M.~Velasco.
\newblock A fractional approach to the {S}turm-{L}iouville problem.
\newblock {\em Central European Journal of Physics}, 11(-):1246--1254, 2013.

\bibitem{marichev}
S.~G. Samko, A.~A. Kilbas, and O.~I. Marichev.
\newblock {\em Fractional integrals and derivatives}.
\newblock Gordon and Breach Science Publishers, Yverdon, 1993.
\newblock Theory and applications, Edited and with a foreword by S. M.
  Nikolski\u{\i}, Translated from the 1987 Russian original, Revised by the
  authors.

\bibitem{Selvadurai}
A.~P.~S. Selvadurai.
\newblock {\em Partial differential equations in mechanics. 1}.
\newblock Springer-Verlag, Berlin, 2000.
\newblock Fundamentals, Laplace's equation, diffusion equation, wave equation.

\bibitem{Ste}
M.~C. Steinbach.
\newblock On {PDE} solution in transient optimization of gas networks.
\newblock {\em J. Comput. Appl. Math.}, 203(2):345--361, 2007.

\bibitem{Vonbelow}
J.~Von~Below.
\newblock Sturm-{L}iouville eigenvalue problems on networks.
\newblock {\em Math. Methods Appl. Sci.}, 10(-):383–395, 1988.

\bibitem{Zayer}
M.~Zayernouri and G.~E. Karniadakis.
\newblock Fractional {S}turm-{L}iouville eigen-problems: theory and numerical
  approximation.
\newblock {\em J. Comput. Phys.}, 252:495--517, 2013.

\bibitem{Zettl}
A.~Zettl.
\newblock {\em Sturm-{L}iouville theory}, volume 121 of {\em Mathematical
  Surveys and Monographs}.
\newblock American Mathematical Society, Providence, RI, 2005.

\end{thebibliography}

\end{document}